\RequirePackage{fix-cm, amsmath}
\documentclass{amsart}     
\usepackage{graphicx}
\usepackage{latexsym}
\usepackage{amssymb}
\usepackage{graphics}
\usepackage{url}
\usepackage[vcentermath]{youngtab}

\usepackage{booktabs}  
\usepackage{ytableau,genyoungtabtikz,graphicx,wrapfig}
\newtheorem{theorem}{Theorem}[section]
\newtheorem*{theorem*}{Theorem}
\newtheorem{corollary}[theorem]{Corollary}

\newtheorem{proposition}[theorem]{Proposition}

\newtheorem{lemma}[theorem]{Lemma}
\newtheorem*{lemma*}{Lemma}
\theoremstyle{definition}
\newtheorem*{definition}{Definition}

\newtheorem{example}[theorem]{Example}

\newcommand{\N}{\mathbf{N}}
\newcommand{\Z}{\mathbf{Z}}
\newcommand{\Q}{\mathbf{Q}}

\newcommand{\F}{\mathbf{F}}

\renewcommand{\epsilon}{\varepsilon}

\DeclareMathOperator{\Sym}{Sym}

\DeclareMathOperator{\Hom}{Hom}

\DeclareMathOperator{\sgn}{sgn}
\DeclareMathOperator{\Inf}{Inf}
\DeclareMathOperator{\Lin}{Lin}

\newcommand{\Ind}{\big\uparrow}
\newcommand{\Res}{\big\downarrow}
\newcommand{\ind}{\!\uparrow}
\newcommand{\res}{\!\downarrow}

\newcommand{\size}[1]{\vert #1 \vert}
\newcommand{\Irr}{\text{\normalfont Irr}}
\newcommand{\supp}{\text{\normalfont supp}}
\newcommand{\tup}[1]{\;\left[#1,\bar{#1}\right]\;}
\numberwithin{equation}{section}
\begin{document}

\title[Vertices and decomposition numbers]{Vertices of modules and decomposition numbers of $C_2\wr S_n$}

\author{Jasdeep Kochhar}


\address{
              Department of Mathematics, Royal Holloway, University of London, Egham, TW20 0EX, United Kingdom}
              \email{Jasdeep.Kochhar.2015@rhul.ac.uk}           

\date{\today}

\begin{abstract}
Given $n \in \mathbf{N},$ we consider the imprimitive wreath product $C_2 \wr S_n.$ 
We study the structure of the $p$-modular reduction of modules whose ordinary characters form an involution model of $C_2 \wr S_n,$ where $p$ is an odd prime.
We describe the vertices of these modules, and we use this description of the vertices to determine certain decomposition numbers of $C_2 \wr S_n.$
\end{abstract}
\maketitle

\section{Introduction}\label{sec: intro}
Given $n \in \N,$ let $S_n$ denote the symmetric group on $n$ letters.
The imprimitive wreath product $C_2 \wr S_n$ can be thought of in various ways, such as the Weyl group of type $B_n,$ or as the symmetry group of the $n$-dimensional hypercube. 
In this paper we use the characterisation of $C_2 \wr S_n$ as the centraliser of the permutation
$$(1\ \bar{1})(2\ \bar{2})\ldots(n\ \bar{n})$$
in the symmetric group on the set $\{1,2,\ldots,n,\bar{1},\bar{2},\ldots,\bar{n}\}.$
Explicitly, $C_2 \wr S_n$ is the subgroup of $S_{2n}$ generated by the set
$$\{(1\ \bar{1}),(1\ 2)(\bar{1}\ \bar{2}),(1\ 2\ldots n)(\bar{1}\ \bar{2}\ldots \bar{n})\}.$$
We study certain representations of $FC_2 \wr S_n,$ where $F$ is a field of odd prime characteristic $p.$
We remark that there are various connections between the representations of $C_2 \wr S_n$ and $S_n.$
For instance, the simple $FC_2 \wr S_n$-modules can be constructed using the simple $FS_n$-modules (see \cite{AMP} and \cite[\textsection 4.3]{JamesKerber}), where $F$ is any field. 
For detailed accounts on the representation theory of the symmetric group, see \cite{JamesKerber} and \cite{James}.

In the case of our attention, there are notable similarities between the block structures of the algebras $FC_2 \wr S_n$ and $FS_n.$
In Proposition \ref{prop: nakayamconj}, we give a complete description of the blocks of $FC_2 \wr S_n.$
We also describe the Brauer correspondence between the blocks of $FC_2 \wr S_n$ and the blocks of $FN_{C_2 \wr S_n}(R),$ where $R$ is a particular $p$-subgroup of $C_2 \wr S_n.$

We also note that if $k$ is a field of characteristic 2, then the subgroup 
$$B_n := \langle (1\ \bar{1}),(2\ \bar{2}),\ldots, (n\ \bar{n})\rangle,$$
is contained in the kernel of all irreducible representations of $kC_2 \wr S_n.$
In this case, the irreducible representations of $kC_2 \wr S_n$ are inflations (see \textsection \ref{sec: main}) of the irreducible representations of $kS_n.$ 
\subsection{Main results}\label{sec: main}
We introduce the background required to state our first main theorem. 

Define $T_n$ to be the subgroup of $C_2 \wr S_n$ such that
$$T_n=\langle(1\ 2)(\overline{1}\ \overline{2}),(1\ 2\ldots n)(\overline{1}\ \overline{2}\ldots \overline{n})\rangle.$$
We refer to $T_n$ as the \emph{top group}.
It follows that $C_2 \wr S_n = B_n \rtimes T_n.$

Consider $\sigma \in \Sym(\{1,2,\ldots,n\}).$ 
Define $\overline{\sigma} \in \Sym(\{\overline{1},\ldots,\overline{n}\})$ to be such that $\overline{\sigma}(\overline{x}) = \overline{\sigma(x)}.$
Given $H \le \Sym(\{1,2,\ldots,n\}),$ we define $\xi(H)$ to be the subgroup of the top group consisting precisely of the permutations $\sigma\overline{\sigma}$ such that $\sigma \in H.$

Let $\hat{h}$ be the image of $h \in C_2 \wr S_n$ under the natural surjection 
$C_2 \wr S_n \twoheadrightarrow S_n.$ 
Write $\Inf_{S_n}^{C_2 \wr S_n} \rho$ for the representation of $FC_2\wr S_n$
such that $(\Inf_{S_n}^{C_2 \wr S_n} \rho)(h) = \rho(\hat{h})$ for $h \in  C_2 \wr S_n,$
where $\rho$ is a representation of $FS_n.$ 
Note that $FB_n$ is in the kernel of $\Inf_{S_n}^{C_2 \wr S_n} \rho.$
If $M$ is a module corresponding to $\rho,$ we write $\Inf_{S_n}^{C_2 \wr S_n} M$ 
for the module corresponding to $\Inf_{S_n}^{C_2 \wr S_n} \rho.$ 

Given $a \in \N,$ define $f_a \in C_2 \wr S_{2a}$ to be the permutation equal to
$$(1\ a+1)(2\ a+2)\ldots(a\ 2a)(\,\overline{1}\ \overline{a+1}\,)(\,\overline{2}\ \overline{a+2}\,)\ldots(\,\overline{a}\ \overline{2a}\,),$$
and let $V_a$ be the centraliser of $f_a$ in $C_2 \wr S_{2a}.$
Therefore $V_a$ is equal to
$$\langle (\,1\ \overline{1}\,)(\,a+1\ \overline{a+1}\,),(\,2\ \overline{2}\,)(\,a+2\ \overline{a+2}\,),\ldots,(\,a\ \overline{a}\,)(\,2a\ \overline{2a}\,)\rangle \rtimes \xi(\,S_2\wr S_a\,).$$
We also define $V_\lambda$ to be the subgroup of $V_a$ equal to 
$$\langle (\,1\ \overline{1}\,)(\,a+1\ \overline{a+1}\,),(\,2\ \overline{2}\,)(\,a+2\ \overline{a+2}\,),\ldots,(\,a\ \overline{a}\,)(\,2a\ \overline{2a}\,)\rangle \rtimes \xi(\,S_2\wr S_\lambda\,),$$
where $\lambda$ is a partition of $a,$ and $S_\lambda$ is the corresponding Young subgroup of $S_a.$ 

Let $N$ denote the non-trivial simple $FC_2$-module. 
Given $c \in \N,$ we define $\widetilde{N}^{\otimes c}$ to be the one-dimensional $FC_2 \wr S_c$-module on which $T_c$ acts trivially and the permutation $(\,i\ \bar{i}\,)$ acts by negative sign for all $1 \le i \le c.$ 
Note that $\widetilde{N}^{\otimes c}$ agrees with the notation in \cite[\textsection 4.3]{JamesKerber}.
Also define $\sgn_{S_c}$ to be the sign module of $FS_c.$
Given $(a,b,c) \in \N^3_0$ such that $2a + b + c = n,$ we define the module
$$M_{(2a,b,c)} = (F\Ind_{V_a}^{C_2 \wr S_{2a}} \boxtimes \Inf_{S_b}^{C_2 \wr S_b} \sgn_{S_b} 
\boxtimes (\widetilde{N}^{\otimes c}\otimes\Inf_{S_c}^{C_2 \wr S_c}\sgn_{S_c}))\Ind_{C_2 \wr S_{(2a,b,c)}}^{C_2 \wr S_n},$$
where $\Ind$ denotes the induction of modules. 
Our first main theorem characterises the vertices of the indecomposable summands of $M_{(2a,b,c)}.$ 
To state this result, we also require the following notation.
Given $r \in \N$ such that $rp \le n,$ define
$$T'_r := \{(\lambda,t,u) : \lambda \in \Lambda(2,s), 2s+t+u = r \mbox{ and } sp \le a, tp \le b, up \le c\},$$
where $\Lambda(2,s)$ denotes the set of all compositions of $s$ in at most 2 parts. 

\begin{theorem}\label{thm: main}
Let $(a,b,c) \in \N^3_0$ be such that $2a + b + c = n,$ and let $U$ be a non-projective indecomposable summand of $M_{(2a,b,c)}$. 
Then $U$ has a vertex equal to a Sylow $p$-subgroup of 
$$V_{p\lambda} \times C_2\wr S_{tp} \times C_2 \wr S_{up},$$
for some $r \in \N,$ where $rp \le n,$ and $(\lambda,t,u) \in T'_r.$
\end{theorem}
A key motivation for Theorem \ref{thm: main} is \cite[Theorem 1.2]{GW},
which describes the vertices of the indecomposable summands of the $FS_{2m+k}$-module
\[H^{(2m;k)}:= (F_{S_2\wr S_m} \boxtimes \sgn_{S_k})\Ind_{S_2\wr S_m \times S_k}^{S_{2m+k}}.\]
This description is used to determine certain decomposition numbers of the symmetric group in \cite[Theorem 1.1]{GW}.
An essential part of the proof of \cite[Theorem 1.1]{GW} is that the ordinary character of $H^{(2m;k)}$ is multiplicity free.
By Propositions 1 and 2 in \cite{baddeley1991models}, the ordinary character of the module 
$M_{(2a,b,c)}$ is also multiplicity free.
Also observe that both of $S_2 \wr S_m \times S_k$ and $V_a \times C_2 \wr S_b \times C_2 \wr S_c$ 
are centralisers of certain involutions in $S_{2m+k}$ and $C_2 \wr S_n,$ respectively.
These similarities between $H^{(2m;k)}$ and $M_{(2a,b,c)}$ arise as these two modules are examples of a deeper phenomenon that is occurring. 
In order to describe this, we require the following definition.  
We say that a finite group $G$ has an {\it involution model} if there exists a set of elements
$\{e_1,e_2,\ldots,e_t\} \subseteq G,$ such that $e_i^2 = 1$ for all $i,$ and for each $e_i$ there exists a linear character $\psi_i$ of $C_G(e_i)$
such that 
$$\sum_{i=1}^t \psi_i^G = \sum_{\psi \in \Irr(G)} \psi.$$
The main result of \cite{inglis1990explicit} is that the sum of the ordinary characters of the modules $H^{(2m;k)}$ is an involution model for $S_{2m+k}.$
In \cite{baddeley1991models}, Baddeley proves that if the group $G$ has an involution model, then the imprimitive wreath product $G \wr S_n$ also has an involution model.
The module $M_{(2a,b,c)}$ is part of the involution model constructed by Baddeley when $G = C_2.$ 
Using similar ideas as in the proof of \cite[Theorem 1.1]{GW}, we use Theorem \ref{thm: main} in this paper to understand particular decomposition numbers of $C_2 \wr S_n.$

In order to state our result on decomposition numbers, we require an understanding of the simple $\Q C_2 \wr S_n$-modules.
Let $\mathcal P^2(n)$ be the set of pairs of partitions $(\lambda,\mu)$ such that $\size{\lambda} + \size{\mu} = n.$ 
Given $(\lambda,\mu) \in \mathcal P^2(n),$ we define 
$$S^{(\lambda,\mu)} = (\Inf_{S_{\size{\lambda}}}^{C_2 \wr S_{\size{\lambda}}}S^\lambda \boxtimes 
\widetilde{N}^{\otimes \size{\mu}}\otimes\Inf_{S_{\size{\mu}}}^{C_2 \wr S_{\size{\mu}}}S^\mu)\Ind_{C_2 \wr S_{(\size{\lambda},\size{\mu})}}^{C_2 \wr S_n},$$
where $S^{\lambda}$ is the usual Specht module labelled by $\lambda.$
We refer to a module of this form as a {\it hyperoctahedral Specht module}.
The set of all hyperoctahedral Specht modules $S^{(\lambda,\mu)}$ such that $(\lambda,\mu) \in \mathcal P^2(n)$ is a complete set of simple $\Q C_2 \wr S_n$-modules.  

Recall that $p$ is a prime number not equal to 2.
Let $(\nu,\widetilde{\nu}) \in \mathcal P^2(n)$ be such that $\nu$ and $\widetilde{\nu}$ are $p$-regular. 
The module $S^{(\nu,\widetilde{\nu})}$ has a unique maximal submodule, and we write $D^{(\nu,\widetilde{\nu})}$ for the quotient of $S^{(\nu,\widetilde{\nu})}$ by this maximal submodule. 
We have that $D^{(\nu,\widetilde{\nu})}$ is a simple $FC_2 \wr S_n$-module; furthermore, every simple $FC_2 \wr S_n$-module
is of this form (see Proposition \ref{prop: simples}). 
The decomposition number $d_{\lambda\nu,\mu\widetilde{\nu}}$ is defined to be the  number of composition factors of
$S^{(\lambda,\mu)}$ isomorphic to $D^{(\nu,\widetilde{\nu})}.$

We now introduce one more piece of notation required to state our second main theorem.
Given a $p$-core partition $\gamma$ and given $b \in \N_0,$ let $w_b(\gamma)$ be the minimum number of $p$-hooks
such that when added to $\gamma,$ we obtain a partition with exactly $b$ odd parts.
Let $\mathcal E_b(\gamma)$ be the set of all partitions of $\size{\gamma} + w_b(\gamma)p$ obtained in this way.

\begin{theorem}\label{thm: main1}
Let $\gamma$ and $\delta$ be $p$-core partitions, and let $b,c \in \N_0.$
If $b \ge p$ (resp. $c \ge p$), suppose that $w_{b-p}(\gamma) \neq w_b(\gamma) - 1$ (resp.~$w_{c-p}(\delta) \neq w_c(\delta)-1).$
Then there exists a set partition of $\mathcal{E}_b(\gamma) \times \mathcal{E}_c(\delta)$, say $\Lambda_1,\ldots,\Lambda_t,$ such that each $\Lambda_i$ has a unique pair $(\nu_i,\widetilde{\nu_i})$ with $\nu_i$ and $\widetilde{\nu_i}$ both maximal in the dominance orders on $\mathcal{E}_b(\gamma)$ and $\mathcal{E}_c(\delta),$ respectively. 

Moreover, $\nu_i$ and $\widetilde{\nu_i}$ are $p$-regular for each $i,$ and the column of the decomposition matrix of $C_2 \wr S_n$ labelled by $(\nu_i,\widetilde{\nu_i})$ has ones in the rows labelled by the pairs in $\Lambda_i$, and zeros in all other rows. 
\end{theorem}
We remark that Theorem \ref{thm: main1} follows by combining \cite[Theorem 1.1]{GW} and the Morita equivalence between $FC_2 \wr S_n$ and $\bigoplus_{i=0}^n FS_{(i,n-i)}$ given by Proposition \ref{prop: nakayamconj} in this paper.
However, our characterisation of the vertices in Theorem \ref{thm: main} does not follow from this Morita equivalence. 
Example \ref{ex: main} makes this explicit by describing the vertices of the non-projective indecomposable summands of $M_{(54,0,0)}$ over a field of characteristic 3.

\subsection{Outline}\label{sec: outline}
We now provide an outline of the paper.
Our approach to proving Theorem \ref{thm: main} is by using results on the Brauer morphism
for $p$-permutation modules (see \cite{broue1985scott}). 
We present the required background on the Brauer morphism in \textsection \ref{sec: Broue}.
In \textsection \ref{sec: top_group} we construct particular subgroups of $C_2 \wr S_n$ that will 
be used in the proof of Theorem \ref{thm: main}.
We also require a description of the conjugacy classes in $C_2 \wr S_n,$ which we give in \textsection \ref{sec: conj}.

We collect the required background for Theorem \ref{thm: main1} in \textsection \ref{sec: hypspecht} and \textsection \ref{sec: blockC2wrSn}.
In \textsection \ref{sec: hypspecht} we describe a basis of $S^{(\lambda,\mu)}.$
In order to apply Theorem \ref{thm: main} to Theorem \ref{thm: main1}, we require results on the blocks
of the group algebra $FC_2 \wr S_n.$ 
We state these results on blocks in \textsection\ref{sec: blockC2wrSn} as a corollary of Theorem \ref{thm: Morita}. 
As a further corollary of Theorem \ref{thm: Morita}, we also describe the simple $FC_2 \wr S_n$-modules in \S\ref{sec: blockC2wrSn}. 
 
In \textsection \ref{sec: Mabc} we give an explicit combinatorial basis for the module $M_{(2a,b,c)},$ specifically in \textsection \ref{sec: basis}. 
The basis that we describe is generally not a permutation basis for $M_{(2a,b,c)},$ and so is in general not a $p$-permutation basis for an arbitrary $p$-subgroup of $C_2 \wr S_n.$
In \textsection \ref{sec: pbasis}, we show how the basis given in \textsection \ref{sec: basis} can be used to construct a $p$-permutation basis of $M_{(2a,b,c)}$  with respect to a given $p$-subgroup of $C_2 \wr S_n.$
We will use this basis and results on the Brauer morphism to prove Theorem \ref{thm: main}.

In \textsection \ref{sec: proof} we prove Theorem \ref{thm: main}. 
The proof is technical in areas, and so it is broken down into three steps. 
We first consider the Brauer correspondent of $M_{(2a,b,c)}$ with respect to a particular cyclic group of order $p$ in $C_2 \wr S_n,$ denoted $R_r,$ where $rp \le n.$
We decompose $M_{(2a,b,c)}$ as a direct sum of indecomposable $FN_{C_2 \wr S_n}(R_r)$-modules, denoted $N_{(\lambda,t,u)},$ using Clifford theory arguments. 
We see that each summand $N_{(\lambda,t,u)}$ has a vertex containing the subgroup $R_{\omega^\star}$
(defined in the first step of the proof).
In the second step, we therefore consider the module $N_{(\lambda,t,u)}(R_{\omega^\star}).$
We show that $N_{(\lambda,t,u)}(R_{\omega^\star})$ is an indecomposable $FN_{C_2 \wr S_n}(R_{\omega^\star})$-module, and we determine its vertex. 
In the third step, we use the description of the vertices of $N_{(\lambda,t,u)}(R_{\omega^\star})$ to complete the proof of Theorem \ref{thm: main}. 

In \textsection \ref{sec: decomp} we prove Theorem \ref{thm: main1}.
We do this by showing that every summand of $M_{(2a,b,c)}$ in the block $B((\gamma,w_b(\gamma)),(\delta,w_c(\delta)))$
is projective. 
We can lift summands of $M_{(2a,b,c)}$ from $\F_p$ to $\Z_p$ using a result of Scott. 
We are then able to understand the ordinary characters of these lifted summands using Brauer reciprocity for projective modules. 
These results of Scott and Brauer are essential tools in using the local information on vertices to understand the situation in the global case. 

\section{Background}
Recall that $F$ is assumed to be a field of characteristic $p,$ where $p$ is an odd prime.
We remark that the results in \textsection \ref{sec: Broue} also hold over a field of characteristic 2.
\subsection{The Brauer morphism}\label{sec: Broue}
Throughout this section let $G$ be a finite group. 
We recount the theory of the Brauer morphism that will be used in this paper.
We also state the results required that relate the Brauer morphism to the vertices of modules and to the blocks of $FG.$
For details on the theory of vertices and the theory of blocks, we refer the reader to \cite{Alperin}.

Let $H \le G$, and let $M$ be an $FG$-module. Define $M^H$ to be the set of vectors in $M$ fixed by every element in $H.$ Given $L \le H \le G$, we define the map
$$\begin{array}{rcl}
\text{Tr}_L^H: M^L & \rightarrow & M^H \\
x& \mapsto &\sum xg,
\end{array}$$
where the sum runs over a transversal of the cosets of $L$ in $H$. 

When $P$ is a $p$-subgroup of $G$, we define  
$$M(P) = M^P/\sum_{Q < P} \text{Tr}_Q^P M^P.$$
This is an $FN_G(P)$-module, on which $P$ acts trivially. 
The natural quotient map $M^P \mapsto M(P),$ is known as the {\it Brauer morphism} and is an $FN_G(P)$-module homomorphism.

Let $U$ be an indecomposable $FG$-module. A vertex of $U$ is a subgroup $Q$ of $G$, minimal such that $U$ is a direct summand of $U\Ind^G_Q$. It is known (see \cite[Section 9, Theorem 4]{Alperin}) that $U$ has a vertex, and that any vertex of $U$ is a $p$-group. Furthermore, any two vertices of $U$ are conjugate in $G$.

The Brauer morphism can be used to obtain information about the vertices of a certain class of modules known as $p$-permutation modules. The module $M$ is a {\it $p$-permutation module} if for all $p$-subgroups of $G$, there exists an $F$-basis of $M$ that is permuted by $P$. If $\mathcal B$ is such a basis, then we say that $\mathcal B$ is a {\it $p$-permutation basis} of $M$ with respect to $P$. The following lemma from \cite{broue1985scott} gives an another characterisation of $p$-permutation modules. 
\begin{lemma}[{\cite[0.4]{broue1985scott}}]
The module $M$ is a $p$-permutation module if and only if there exists a subgroup $H$ of $G$ such that $M$ is a summand of $F\Ind_H^G.$ 
\end{lemma}
We now assume that $M$ is a $p$-permutation module, and that $P$ is a $p$-subgroup of $G.$
The following lemmas show how the Brauer morphism can be used to determine the vertices of a $p$-permutation module.
\begin{lemma}[{\cite[3.2(1)]{broue1985scott}}]\label{lem: brauer_vertex}
Let $M$ be an indecomposable $p$-permutation $FG$-module.
Then $M$ has a vertex equal to $P$ if and only if $P$ is a maximal $p$-subgroup of $G$ such that $M(P) \neq 0.$
\end{lemma}
\begin{lemma}[{\cite[1.1(3)]{broue1985scott}}]\label{lem: brauer_fixed}
Let $\mathcal{B}$ be a $p$-permutation basis of $M$ with respect to $P,$ and let $\mathcal B^P$ be the set of points in $\mathcal B$ that are fixed by $P$.
Then $\mathcal B^P$ is a basis of $M(P).$
\end{lemma}
It follows that $Q$ is a vertex of $M$ if there exists a vector in a $p$-permutation basis of $M$ (with respect to $Q$) that is fixed by $Q,$ and $Q$ is maximal with this property. 
We also require the following lemma, which will be crucial in the proof of Theorem \ref{thm: main}.
\begin{lemma}[{\cite[Lemma 4.7]{GW}}]\label{lem: broue_norm}
Let $R$ be a normal subgroup of $P,$ and let $K = N_G(R).$
Then $M(R)$ is a $p$-permutation $FK$-module. Moreover,
$M(P) \cong M(R)(P),$ where the isomorphism is of $FN_K(P)$-modules.
\end{lemma}

We now describe how the Brauer morphism can be used to determine information about blocks of a group algebra.
Given $H \le G,$ let $B$ a block of $G,$ and let $b$ be a block of $H$. 
We say that the block $B$ {\it corresponds to} $b$ (see \cite[\textsection 4]{Alperin}) if $b$ is a summand of $B_{H \times H},$ and $B$ is the unique block of $G$ with this property. 
In this case, we write $b^G = B.$ 

The following lemma relates the theory of blocks to the Brauer morphism. 
\begin{lemma}[{\cite[Lemma 7.4]{wildon2010vertices}}]\label{lem: broue_covering}
Suppose that $M$ lies in the block $B$ of $G$. 
If $M(P)$ has a summand in the block $b$ of $N_G(P),$ then $b^G = B.$
\end{lemma}
\subsection{Subgroups of $C_2 \wr S_n$}\label{sec: top_group}
Recall that given an element $h \in C_2 \wr S_n,$ we write $\widehat{h}$ for the image of $h$ under the natural surjection $C_2 \wr S_n \twoheadrightarrow S_n.$ 
Given $Q \le C_2 \wr S_n$, we define $\widehat{Q} = \{\widehat{h} : h \in Q\}.$

Also given $X \subset \{1,2,\ldots,n\},$ we write $C_2 \wr S_X$ for the subgroup of $C_2\wr S_n$ generated by the set
$$\{(x\ \overline{x}): x \in X\} \cup \{(x\ y)(\overline{x}\ \overline{y}) : x,y \in X, x\neq y\}.$$

We now consider $p$-subgroups of $C_2 \wr S_n.$
As $C_2 \wr S_n$ has a Sylow $p$-subgroup contained in the top group $T_n,$ any $p$-subgroup of $C_2 \wr S_n$ has a conjugate in the top group. 
\subsection{Conjugacy in the hyperoctahedral group}\label{sec: conj}
Given $i\in \{1,2,\ldots, n\}$, we define $\overline{\overline{i}} = i.$ 
Given $g \in C_2 \wr S_n$, we say that $g$ is a {\it positive $r$-cycle} if 
$$g = (a_1,a_2,\ldots,a_r)(\overline{a_1},\overline{a_2},\ldots,\overline{a_r}),$$
and that $g$ is a {\it negative $r$-cycle} if 
$$g = (a_1,a_2,\ldots,a_r,\overline{a_1},\overline{a_2},\ldots,\overline{a_r}),$$
where $a_1,\ldots,a_r \in \{1,\overline{1},\ldots,{n},\overline{n}\}.$
\begin{example}
Let $n = 1.$
The identity permutation $(1)(\bar{1})$ is a positive 1-cycle, and the permutation $(1\ \overline{1})$ is a negative 1-cycle.
\end{example}
Every element of $C_2 \wr S_n$ can be expressed uniquely, up to the order of the factors, as a product of disjoint positive and negative cycles. 
The number of positive (resp. negative) $r$-cycles of $g \in C_2 \wr S_n$ is denoted by $p_r$ (resp. $n_r$), and we say that $g$ has {\it cycle type} $((p_r),(n_r))_{1\le r \le n}.$
We now have the following lemma.
\begin{lemma}[{\cite[\textsection 4.2]{JamesKerber}}]\label{lem: conj}
Let $g, h \in C_2 \wr S_n.$ Then $g$ and $h$ are conjugate in $C_2 \wr S_n$ if and only if they have the same cycle type. 
\end{lemma}
Furthermore, the centraliser of an element with cycle type $((p_r),(n_r))_{1\le r \le {n}}$ has order equal to 
$$\prod_{r = 1}^n(2r)^{p_r+n_r}(p_r!)(n_r!).$$
\subsection{Hyperoctahedral Specht modules}\label{sec: hypspecht}
Given $x\in\{1,2,\ldots, n\},$ we define $[x,\overline{x}]$ to be the image of $(x,\overline{x})$ in the quotient of the $FC_2 \wr S_n$-permutation module $ F[\{1,\ldots,n,\overline{1},\ldots,\overline{n}\}]$ by the submodule generated by
$$\{(x,\overline{x}) + (\overline{x},x) : 1 \le x \le n\}.$$
Therefore the $F$-span of $[x,\overline{x}]$ is isomorphic to $N$ as an  $F[\Sym(\{x,\overline{x}\})]$-module.

Given $(\lambda,\mu) \in \mathcal P^2(n),$ let $t$ be the disjoint union of a $\lambda$-tableau and a $\mu$-tableau, such that 
\begin{enumerate}
\item the $\lambda$-tableau has entries $\{x,\overline{x}\}$, and the $\mu$-tableau has entries $[y,\overline{y}]$
\item the set $\{x,\overline{x}\}$ is an entry of the $\lambda$-tableau if and only if $[x,\overline{x}]$ is not an entry of the $\mu$-tableau, for all $1 \le x \le n.$
\end{enumerate}

In this case, we say that $t$ is a $(\lambda,\mu)$-tableau.
We write $t^+$ for the $\lambda$-tableau, and $t^-$ for the $\mu$-tableau.
\begin{example}
The following is a $((3),(3,1))$-tableau. 
$$\ytableausetup
{mathmode, boxsize=1.6em}
\begin{ytableau}
\none&\none&\none&\scriptstyle {\tup{4}}& \scriptstyle {\tup{5}}&\scriptstyle {\tup{6}}\\
\none&\none&\none&\scriptstyle {\tup{7}}\\
\scriptstyle {\{1,\bar{1}\}}&\scriptstyle {\{2,\bar{2}\}}&\scriptstyle {\{3,\bar{3}\}}\\
\end{ytableau}
$$

\end{example}

Given a $(\lambda,\mu)$-tableau $t,$ let $R(t)$ (resp.~$C(t)$) be the subgroup of $T_n$ consisting of all permutations that setwise fix the entries in each row (resp.~column) of $t$. 
We define an equivalence relation $\backsim$ on the set of $(\lambda,\mu)$-tableaux by $t \backsim u$ if and only if there exists $\pi \in R(t)$ such that $u =  t\pi.$
The {\it $(\lambda,\mu)$-tabloid} $\{t\}$ is the equivalence class of $t$. 
We define the \emph{$(\lambda,\mu)$-polytabloid}  $e_t$ by
\[ e_t = \sum_{\sigma\in C(t)} \text{sgn}(\sigma)\{t\}\sigma. \]
The hyperoctahedral Specht module $S^{(\lambda,\mu)}$ (defined in \textsection \ref{sec: outline}) is a cyclic module, generated by any $(\lambda,\mu)$-polytabloid $e_t$.
In order to describe a basis of $S^{(\lambda,\mu)}$, we order the sets $\{x,\overline{x}\}$ by setting $\{x,\overline{x}\} \le \{y,\overline{y}\}$
if and only if $x \le y.$
We also define an ordering on the set of $[x,\overline{x}]$ in the same way. 
We say that $t$ is {\it standard} if both $t^+$ and $t^-$ are standard tableaux with respect to the orders just defined. 
The module $S^{(\lambda,\mu)}$ has a basis given by the set of all polytabloids $e_t$ such that $t$ is a standard $(\lambda,\mu)$-tableau.
\subsection{The simple modules and the blocks of $FC_2 \wr S_n$}\label{sec: blockC2wrSn}
It is well-known that the blocks of $S_n$ are labelled by pairs $(\gamma,v)$ such that $\gamma$ is a $p$-core partition, and
$\vert \gamma \vert + vp = n.$
Moreover, the $FS_n$-module $S^{\lambda}$ lies in the block labelled by the $p$-core of $\lambda.$
This result is known as Nakayama's conjecture, and was first proved by Brauer and Robinson in \cite{brauernakay} and \cite{robinnakay}.
The main result in this section is a complete description of the blocks of $FC_2 \wr S_n,$ which we give in Proposition \ref{prop: nakayamconj}. 
We actually prove the stronger Theorem \ref{thm: Morita} below, from which Proposition \ref{prop: nakayamconj} follows.
We prove Theorem \ref{thm: Morita} as it is also used in this section to describe the simple $FC_2 \wr S_n$-modules, and in \textsection \ref{sec: decomp} to determine the blocks of $N_{C_2 \wr S_n}(R_r),$ 
where $R_r$ is defined in \textsection \ref{sec: proof}.

We now give the required preliminaries for Theorem \ref{thm: Morita}.
Assume that $G = C_2^a \rtimes H,$ where $a \in \N.$ 
There is an action of $G$ on $\Lin(C_2^a)$ given by conjugation, and we have the following lemma. 
\begin{lemma}\label{lem: conj_block}
The $G$-conjugacy classes of $\Lin(C_2^a)$ are labelled by pairs $(a_1,a_2) \in \N_0^2$ such that $a_1 + a_2 = a.$ 
\end{lemma}
Given $0 \le i \le a,$ write $\Lin_i(C_2^a)$ for the conjugacy class of $\Lin(C_2^a)$ labelled by $(i,a-i).$ 
Fix $\chi_i \in \Lin_i(C_2^a)$ and define $G_i = C_2^a\rtimes H_i,$ where $H_i$ is the stabiliser of $\chi_i$ in $H.$  
Given an $FG$-module $V$ and $\chi \in \Lin(C_2^a),$ let 
\[V^\chi = \{v \in V : vg = \chi(g)v \mbox{ for all } g \in C_2^a\}.\]
For $g \in G,$ we have that $V^{\chi}g = V^{\chi g},$ and so $V^{\chi_i}$ is an $FG_i$-module.
Furthermore, $V(i):= \bigoplus_{\chi \in \Lin_i(C_2^a)} V^\chi$ is an $FG$-module.
Then
\begin{equation}\label{eq: weightspacedecomp}
V = \bigoplus_{i=0}^n V(i),
\end{equation}
as a direct sum of $FG$-modules.
We say that $V$ \emph{belongs to} $i$ if $V = V(i)$ for some $i.$ 
Clearly every indecomposable $FG$-module belongs to $i$ for some $i.$  

Let $\theta \in \Hom_{FG}(U,V).$ 
By considering the action of $C_2^a,$ we see that $\theta(U^{\chi}) \subseteq V^{\chi}.$
Therefore $\Hom_{FG}(U,V) = 0$ if $U$ belongs to $i$ and $V$ belongs to $j$ for $i \neq j.$ 
It follows that the $FG$-modules belonging to $i$ generate a subcategory of the module category 
{\bf mod}$(G).$
We write ${\bf mod}_i(G)$ for this subcategory.
\begin{theorem}\label{thm: Morita}
The rings $FG$ and $\bigoplus_{i=0}^{n} FH_i$ are Morita equivalent. 
\end{theorem}
\begin{proof}
Fix $0 \le i \le a.$ 
Let $M$ be an $FH_i$-module, and write $K_i$ for the one-dimensional $FG_i$-module 
on which $C_2^a$ acts according to $\chi_i$ and $H_i$ acts trivially. 
Define the functor $\mathcal F_i: {\bf mod}(H_i) \rightarrow {\bf mod}_i(G)$ by
\[M \mapsto (K_i \otimes \Inf_{H_i}^{G_i} M)\Ind_{G_i}^G.\]
It is sufficient to prove that $\mathcal F_i$ is an equivalence of categories,
which we do by showing that it is essentially surjective, full, and faithful.

To prove that $\mathcal F_i$ is essentially surjective, it is sufficient to consider the case when $U$ is an indecomposable $FG$-module. 
Therefore $U$ belongs to $i,$ and so by definition
\[U = \bigoplus_{\chi \in \Lin_i(C_2^a)} U^{\chi} \cong U^{\chi_i}\Ind_{G_i}^{G}.\]
where the isomorphism follows from \cite[Section 8, Corollary 3]{Alperin}.
By definition, $U^{\chi_i}$ is such that $C_2^a$ acts according to $\chi_i.$ 
Therefore $U^{\chi_i}$ is isomorphic to the tensor product of $K_i$ and a module on which $C_2^a$ acts trivially.
This is equivalent to writing $U^{\chi_i} \cong K_i \otimes \Inf_{H_i}^{G_i} U',$ where $U'$ is an $FH_i$-module.
This proves that $\mathcal F_i$ is essentially surjective. 

Suppose that $0 \neq \theta \in \Hom_{FG}(U,V),$ where $V$ also belongs to $i.$
Write $\varphi$ for $\theta$ restricted to $U^{\chi_i},$ which we view as an $FG_i$-module homomorphism.
We have that $U$ is generated by $U^{\chi_i},$ and so $\varphi(U^{\chi_i}) \neq 0.$ 
Moreover, let $u \in U$ be such that $u = u'g$ for some $g \in G_i/G$ and $u' \in U^{\chi_i}.$
By the remark preceding this proof, we have $\varphi(U^{\chi_i}) \subseteq V^{\chi_i}.$
Furthermore, by the discussion in the previous paragraph, we have that $U^{\chi_i} \cong U'$ as an $FH_i$-module. 
Writing $\varphi'$ for $\varphi$ viewed as an $FH_i$-module homomorphism, we have
\[\theta(u) = \theta(u'g) = \theta(u')g = \varphi(u')g = \varphi'(u')g.\]
It follows from part (4) of \cite[Section 8, Lemma 6]{Alperin} that $\theta = \mathcal F_i(\varphi'),$ and so $\mathcal F_i$ is full. 
Moreover, $\varphi'$ is determined by the restriction of $\theta$ to $U^{\chi_i},$ and so $\mathcal F_i$ is faithful. 
\end{proof}

\begin{proposition}\label{prop: nakayamconj}
The rings $FC_2 \wr S_n$ and $\bigoplus_{i=0}^n FS_{(i,n-i)}$ are Morita equivalent. 
Moreover, the blocks of $FC_2 \wr S_n$ are labelled by pairs $((\gamma,v),(\delta,w)),$ 
where $\gamma$ and $\delta$ are $p$-core partitions such that $\size{\gamma} + vp + \size{\delta}+wp = n.$
The hyperoctahedral Specht module $S^{(\lambda,\mu)}$ lies in the block labelled by $((\gamma,v),(\delta,w))$ if and only if 
$\lambda$ is a partition of $\size{\gamma}+vp$ with $p$-core $\gamma,$ and $\mu$ is a partition of $\size{\delta}+wp$ with $p$-core $\delta.$
\end{proposition}
\begin{proof}
Given $i \in \{0,1,\ldots,n\},$ let $\chi_i \in \Lin(C_2^n)$ be such that $\chi_i((1\ \overline{1})) = \cdots = \chi_i((i\ \overline{i})) = 1$ 
and $\chi_i((i+1\ \overline{i+1})) = \cdots = \chi_i((n\ \bar{n})) = -1.$ 
In this case $H_i = S_{(i,n-i)}.$
The first statement of the result is now immediate using Theorem \ref{thm: Morita}. 
The remaining statements follow from the definition of $S^{(\lambda,\mu)}$ and Nakayama's conjecture for the symmetric group. 
\end{proof}
We write $B((\gamma,v),(\delta,w))$ for the block labelled by the pair $((\gamma,v),(\delta,w)).$ 

We now describe the simple $FC_2 \wr S_n$-modules $D^{(\nu,\widetilde{\nu})}.$ 
Given a $p$-regular partition $\nu,$ the Specht module $S^{\nu}$ has a unique maximal submodule, 
and the quotient of $S^{\nu}$ by its maximal submodule is denoted $D^{\nu}.$ 
Moreover, $D^{\nu}$ is a simple $FS_{\size{\nu}}$-module, and every simple $FS_{\size{\nu}}$-module is of this form (see \cite[Theorem 11.5]{James}).

Given $(\nu,\widetilde{\nu}) \in \mathcal P^2(n)$ such that $\nu$ and $\widetilde{\nu}$ are $p$-regular, we define 
\[D^{(\nu,\widetilde{\nu})} = (\Inf_{S_{\size{\nu}}}^{C_2 \wr S_{\size{\nu}}} D^{\nu} \boxtimes 
N^{\otimes \size{\widetilde{\nu}}}\otimes \Inf_{S_{\size{\widetilde{\nu}}}}^{C_2 \wr S_{\size{\widetilde{\nu}}}} D^{\widetilde{\nu}})
\Ind_{C_2 \wr S_{(\size{\nu},\size{\widetilde{\nu}}})}^{C_2 \wr S_n}.\]
We remark that this definition of $D^{(\nu,\widetilde{\nu})}$ agrees with that given in \textsection \ref{sec: main}. 
The following proposition follows immediately from the first statement of Proposition \ref{prop: nakayamconj}. 
\begin{proposition}\label{prop: simples}
Let $n \in \N.$
The set
$$\{D^{(\nu,\widetilde{\nu})} : (\nu,\widetilde{\nu}) \in \mathcal P^2(n) \mbox{ and } \nu,\widetilde{\nu} \mbox{ are $p$-regular} \},$$
is a complete set of non-isomorphic simple $FC_2 \wr S_n$-modules. 
\end{proposition}
\section{A $p$-permutation basis of $M_{(2a,b,c)}$}\label{sec: Mabc}
Let $a,b,c\in \N_0$ be such that $n = 2a+b+c.$ 
In this section we explicitly construct the module $M_{(2a,b,c)}.$ 
We also provide a $p$-permutation basis of $M_{(2a,b,c)}$ with respect to an arbitrary $p$-subgroup of $C_2 \wr S_n.$
\subsection{A module isomorphic to $M_{(2a,b,c)}$}\label{sec: basis}
Let $\mathcal C_{(2a,b,c)}$ be the set
$$\left\lbrace\{g, \gamma, \delta\} :\ 
\begin{array}{l}
\mbox{$g \in C_2 \wr S_n$ has cycle type $a$ positive 2-cycles}\\ 
\gamma = (\{i_{a+1},\overline{i_{a+1}}\},\ldots,\{i_{a+b},\overline{i_{a+b}}\})\\
\delta = ([i_{a+b+1},\overline{i_{a+b+1}}],\ldots,[i_n,\overline{i_n}])\\
\supp(g) \cup \{i_{a+1},\overline{i_{a+1}},\ldots,i_n,\overline{i_n}\} = \{1,\overline{1},\ldots,n,\overline{n}\}
\end{array}
\right\rbrace,$$
where $[x,\overline{x}] = -[\overline{x},x]$ as in \textsection \ref{sec: hypspecht}.

Let $v = \{g,\gamma,\delta\} \in \mathcal C_{(2a,b,c)}$ be such that
\begin{align*}
\gamma &= (\{i_{a+1},\overline{i_{a+1}}\},\ldots,\{i_{a+b},\overline{i_{a+b}}\})\\
\delta &= ([i_{a+b+1},\overline{i_{a+b+1}}],\ldots,[i_n,\overline{i_n}]).
\end{align*}
Define
\begin{align*}
\mathcal S(v) &= \supp(g) \cap \{1,2,\ldots,n\}\\
\mathcal T(v) &= \{i_{a+1},\ldots,i_{a+b}\}\\
\mathcal U(v) &= \{i_{a+b+1},\ldots,i_n\}.
\end{align*} 
As $2a+b+c=n,$ these sets are mutually disjoint.
 
There is an action of $h \in C_2 \wr S_n$ on $v$ given by $vh = \{g^h, \gamma h,\delta h\}.$
With $\mathcal D_{(2a,b,c)}$ defined to be $F$-span of the set
$$\{v - \sgn{(\,\widehat{h}\,)}vh : v \in \mathcal C_{(2a,b,c)},h\in  C_2 \wr S_{\mathcal T(v)} \times C_2 \wr S_{\mathcal U(v)}\},$$
we have the following lemma. 
\begin{lemma}\label{lem: submod}
The vector space $F\mathcal D_{(2a,b,c)}$ is an $FC_2\wr S_n$-submodule of $F\mathcal C_{(2a,b,c)}.$  
\end{lemma}
\begin{proof}
We show that $F\mathcal D_{(2a,b,c)}$ is closed under the action of $C_2 \wr S_n.$ Let $h \in C_2 \wr S_n,$ and let 
$v - \sgn(\widehat g)vg \in \mathcal D_{(2a,b,c)},$
where $g \in C_2 \wr S_{\mathcal T(v)} \times C_2 \wr S_{\mathcal U(v)}.$  
With $g' := g^h,$ it follows that
\begin{align*}
(v - \sgn(\widehat g)vg)h &= vh - \sgn(\widehat g)v(gh)\\
&= vh -\sgn(\widehat g)(vh)g'\\
&= vh - \sgn(\widehat g')(vh)g'.
\end{align*}
The third equality follows as $g$ and $g'$ are conjugate in $C_2 \wr S_n.$
By definition of $\mathcal T(v)$, the set $\mathcal T(vh) = \{x\widehat{h} : x \in \mathcal T(v)\},$
and the analogous statement holds for $\mathcal U(v).$ 
The lemma is now proved as $\supp(g') = \{xh : x \in \supp(g)\},$ and so $g' \in C_2 \wr S_{\mathcal T(vh)} \times C_2 \wr S_{\mathcal U(vh)}.$ 
\end{proof}
Let $v = \{g,\gamma,\delta\} +\mathcal D_{(2a,b,c)}$ be such that
\begin{align*}
\gamma &= (\{i_{a+1},\overline{i_{a+1}}\},\ldots,\{i_{a+b},\overline{i_{a+b}}\})\\
\delta &= ([i_{a+b+1},\overline{i_{a+b+1}}],\ldots,[i_n,\overline{i_n}]),
\end{align*}
where $i_{a+1},\ldots,i_n \in \{1,2,\ldots,n\},$ with $i_{a+1} < \cdots < i_{a+b}$ and $i_{a+b+1} < \cdots <  i_n.$
Write $\mathcal B_{(2a,b,c)}$ for the set of all $v + \mathcal D_{(2a,b,c)}$ of this form.
It follows from Lemma \ref{lem: submod} that $\mathcal B_{(2a,b,c)}$ is a basis of $F\mathcal C_{(2a,b,c)}/\mathcal D_{(2a,b,c)}.$
We use this basis in the following lemma to show that the quotient module $F\mathcal C_{(2a,b,c)}/\mathcal D_{(2a,b,c)}$ is isomorphic to $M_{(2a,b,c)}$ as an $FC_2\wr S_n$-module. 
To simplify the notation, we write $(g,\gamma,\delta)$ for $\{g,\gamma,\delta\} + \mathcal D_{(2a,b,c)} \in \mathcal B_{(2a,b,c)}.$
\begin{lemma} 
The $F$-span of $\mathcal B_{(2a,b,c)}$ is isomorphic to $M_{(2a,b,c)}$ as an $FC_2 \wr S_n$-module.
\end{lemma}
\begin{proof}
Recall that $f_a$ is the element equal to 
$$(1\ a+1)(2\ a+2)\ldots(a\ 2a)(\,\overline{1}\ \overline{a+1}\,)(\,\overline{2}\ \overline{a+2}\,)\ldots(\,\overline{a}\ \overline{2a}\,),$$
with centraliser $V_a$ in $C_2\wr S_{2a}.$
It follows that the module $F\Ind_{V_a}^{C_2\wr S_{2a}}$ has a basis indexed by the elements in the conjugacy class of $f_a$ in $C_2 \wr S_{2a}.$ Let 
\begin{align*}
\gamma &= (\{2a+1,\overline{2a+1}\},\ldots,\{2a+b,\overline{2a+b}\})\\
\delta &= ([2a+b+1,\overline{2a+b+1}],\ldots,[n,\overline{n}]),
\end{align*}
and define $S$ to be the $F$-span of $\{(f_a^g,\gamma ,\delta) : g \in C_2 \wr S_{2a}\}.$ Then $S$ is isomorphic, as an $F[C_2\wr (S_{\{1,2,\ldots,2a\}}\times S_{\{2a+1,\ldots,2a+b\}}\times S_{\{2a+b+1,\ldots,n\}})]$-module, to
$$F\Ind_{V_a}^{C_2\wr S_{2a}}\boxtimes(\Inf_{S_b}^{C_2\wr S_b}\sgn_{S_b}) \boxtimes (\widetilde{N}^{\otimes c} \otimes \Inf_{S_c}^{C_2\wr S_c}\sgn_{S_c}),$$
where $N$ is the non-trivial one-dimensional $FC_2$-module. 
Let 
$$w = (h, (\{j_{a+1},\overline{j_{a+1}}\},\ldots,\{j_{a+b},\overline{j_{a+b}}\}),([j_{a+b+1},\overline{j_{a+b+1}}],\ldots,[j_n,\overline{j_n}])),$$ 
be a vector in $\mathcal B_{(2a,b,c)}.$ As the natural action of $C_2 \wr S_n$ on its blocks 
$$\{1,\overline{1}\},\{2,\overline{2}\},\ldots,\{n,\overline{n}\},$$
is transitive, there exists $\sigma \in C_2\wr S_n$ such that $f_a^\sigma = h,$ and $k^\sigma = j_k$ for all $a+1 \le k \le n.$ It follows that $\overline{v}\sigma = \pm \overline{w},$ and so $FS$ generates $F\mathcal C_{(2a,b,c)}/\mathcal D_{(2a,b,c)}.$ 
Recall that $F\mathcal C_{(2a,b,c)}/\mathcal D_{(2a,b,c)}$ has a basis indexed by elements of the form $(f_a^g,\tilde{\gamma},\tilde{\delta}).$
Using the remark following Lemma \ref{lem: conj}, there are 
$$ \dfrac{2^nn!}{4^aa!2^{b+c}(b+c)!}$$
conjugates of $f_a$ in $C_2 \wr S_n.$ 
Given any such conjugate there are 
$$\binom{b+c}{b}$$
ways to choose the support of $\widetilde{\gamma},$ which then determines $\widetilde{\gamma}$ and $\widetilde{\delta}$ completely.
Therefore
\begin{align*}
\dim_F M_{(2a,b,c)} &= \dfrac{2^nn!}{4^aa!2^{b+c}(b+c)!} \times \binom{b+c}{b}\\
&= \dfrac{2^{2a}\times (2a)!}{4^aa!} \times \dfrac{n!}{(2a)!b!c!}\\
&= \dim_F FS \times [C_2 \wr S_n: C_2 \wr S_{(2a, b, c)}].
\end{align*}
The result now follows by the characterisation of induced modules in \cite[Section 8, Corollary 3]{Alperin}.
\end{proof}
Consider now the module $M_{(2a,0,0)}$, which is a permutation module and therefore a $p$-permutation module. 
The modules $M_{(0,b,0)}$ and  $M_{(0,0,c)}$ are one-dimensional modules. 
Therefore the action of any $p$-subgroup of $C_2 \wr S_b$ or $C_2\wr S_c$ on $M_{(0,b,0)}$ or $M_{(0,0,c)}$, respectively, is trivial. 
It follows that both $M_{(0,b,0)}$ and $M_{(0,0,c)}$ are $p$-permutation modules. 
By \cite[Proposition 0.2(2)]{broue1985scott}, the module $M_{(2a,b,c)}$ is therefore a $p$-permutation module. 

\subsection{A $p$-permutation basis of $M_{(2a,b,c)}$}\label{sec: pbasis}
In this section we assume that $Q$ is a $p$-group contained in the top group of $C_2 \wr S_n,$ which we can do by the discussion in \textsection \ref{sec: top_group}.  
Also given $(g,\gamma,\delta) \in \mathcal B_{(2a,b,c)}$ such that 
\begin{align*}
\gamma &= (\{i_{a+1},\overline{i_{a+1}}\},\ldots,\{i_{a+b},\overline{i_{a+b}}\})\\
\delta &= ([i_{a+b+1},\overline{i_{a+b+1}}],\ldots,[i_n,\overline{i_n}]),
\end{align*} 
define $\theta((g,\gamma,\delta)) = (g,\gamma',\delta')$ where
\begin{align*}
\gamma' &= \{\{i_{a+1},\overline{i_{a+1}}\},\ldots,\{i_{a+b},\overline{i_{a+b}}\}\}\\
\delta' &= \{[i_{a+b+1},\overline{i_{a+b+1}}],\ldots,[i_n,\overline{i_n}]\}.
\end{align*} 
\begin{lemma}\label{lem: ppermbasis}
Let $Q$ be a $p$-subgroup of $C_2 \wr S_n.$ 
Then
\begin{enumerate}
\item there is a choice of sign $s_v$ for each $v \in \mathcal B_{(2a,b,c)}$ such that 
$$\{s_v v : v \in \mathcal B_{(2a,b,c)}\}$$
is a $p$-permutation basis of $M_{(2a,b,c)}$ with respect to $Q$
\item the element $v$ is fixed by $Q$ if and only $\theta(v)$ is fixed by $Q.$
In this case, $s_v = 1.$ 
\end{enumerate}
\end{lemma}
\begin{proof}

Let $H_{(2a,b,c)}$ be the set 
$$\left\lbrace \theta(v) : v \in \mathcal B_{(2a,b,c)}\right\rbrace.$$
There exists a natural bijection between $H_{(2a,b,c)}$ and $\mathcal B_{(2a,b,c)},$ and there is a natural action of $Q$ on $H_{(2a,b,c)}.$

Let $\theta(v_1),\theta(v_2),\ldots,\theta(v_l)$ be representatives for the $Q$-orbits on $H_{(2a,b,c)}.$
Given $\theta(v) \in H_{(2a,b,c)},$ there exists a unique $k$ such that $\theta(v_k) = \theta(v)g$ for some $g \in Q.$
Then $vg$ and $v_k$ are equal up to some ordering of the elements in their respective $b$-tuples and $c$-tuples.
Therefore $v_k = s_vvg,$
for some $s_v \in \{-1,+1\}.$ 

Suppose that there exists some other $\widetilde{g} \in Q$ such that $\theta(v_k) = \theta(v)\widetilde{g}.$
Then $\pm v= vg\widetilde{g}^{-1},$
and so the $F$-span of $v$ is a one-dimensional module for the cyclic group generated by $g\widetilde{g}^{-1}.$
The only such module is the trivial module, and so $vg = v\widetilde{g}.$
The sign $s_v$ is therefore well-defined.

In order to prove the first part of the lemma, we need to check that the set
$$\{s_vv : v \in \mathcal B_{(2a,b,c)}\}$$
is a $p$-permutation basis for $M_{(2a,b,c)}$ with respect to $Q.$
Suppose that $h \in Q$ is such that $s_vvh = \pm s_ww,$
for $v$ and $w$ in $\mathcal B_{(2a,b,c)}.$
Then $s_vv$ and $\pm s_ww$ lie in the same $Q$-orbit, and so there exists some $k$ such that $s_vv = v_kg,$
and $\pm s_ww = v_k\widetilde{g}.$
Therefore $v_k gh\widetilde{g}^{-1} = \pm v_k.$ 
Arguing as before shows that the sign on the right hand side is positive, and so the first part of the lemma is proved.

For the second part of the lemma, if 
$$\theta(v):=(g,\{\{i_{a+1},\overline{i_{a+1}}\},\ldots,\{i_{a+b},\overline{i_{a+b}}\}\},\{[i_{a+b+1},\overline{i_{a+b+1}}],\ldots,[i_n,\overline{i_n}]\})$$
is fixed by $Q,$ then $vh = \pm v$ for all $h \in Q.$ 
Therefore the $F$-span of $v$ is a one-dimensional $Q$-module, and so $v$ is fixed by $Q$ as required. 
Moreover, as $\theta(v)$ is its own $Q$-orbit representative, we have that $s_{v} = 1.$ 
\end{proof}
\section{The vertices of the summands of $M_{(2a,b,c)}$}\label{sec: proof}
Let $U$ be a non-projective indecomposable summand of $M_{(2a,b,c)}.$ 
The vertex of $U$ is therefore non-trivial, and so it contains a conjugate of the cyclic group $C_p$ (viewed as a subgroup of $C_2 \wr S_n).$ 
By the discussion in \textsection \ref{sec: top_group}, any copy of $C_p$ in $C_2 \wr S_n$ is conjugate to
$$R_{r} := \langle \sigma_1\sigma_2\ldots \sigma_{r}\rangle,$$
where $\sigma_j := ((j-1)p+1\ \ldots\ jp)(\overline{(j-1)p+1}\ \ldots \ \overline{jp}),$ for some $rp \le n.$ 
It follows that $U(R_{r}) \neq 0,$ and so in the first step of the proof of Theorem \ref{thm: main}, we completely determine the indecomposable summands of $M_{(2a,b,c)}(R_r).$ 
We begin by describing the group $N_{C_2 \wr S_n}(R_r).$
\subsection{The normaliser of $R_r$}\label{sec: normaliser}
It is clear that there is a factorisation
\begin{equation}\label{eqn: normaliser_fact}
N_{C_2 \wr S_{n}}(R_r) = N_{C_2 \wr S_{rp}}(R_r) \times C_2\wr S_{\{rp+1,\ldots,n\}},
\end{equation}
and so it suffices to describe the group $N_{C_2 \wr S_{rp}}(R_r).$

Let $j \in \N$ be such that $j \le r.$ Define 
$$\tau_j= ((j-1)p+1\ \overline{(j-1)p+1})\ldots(jp\ \overline{jp}).$$
The subgroup $\langle \tau_1, \tau_2,\ldots, \tau_r \rangle$ is the full centraliser of $R_r$ in the  subgroup $B_{rp}$ (as defined in \textsection \ref{sec: intro}) of $C_2 \wr S_{rp}.$

Let $i \in \N$ be such that $i < r.$ Define 
$$\rho_i =  ((i-1)p+1\ ip+1)(\overline{(i-1)p+1}\ \overline{ip+1})\ldots (ip \ (i+1)p)(\overline{ip} \ \overline{(i+1)p}).$$
We note that 
$$
\sigma_j^{\rho_i} = \left\lbrace 
\begin{array}{ll} 
\sigma_{j+1} & j = i\\
\sigma_{j-1} & j = i+1\\
 \sigma_j & j \not\in \{i,i+1\}.
\end{array}\right.
$$
Let $x$ be a fixed primitive root modulo $p$. 
Given $i \in \N,$ let $j$ be the unique natural number such that $(j-1)p < i \le jp.$ 
We define $z_r \in C_2 \wr S_{rp}$ to be the permutation such that $z_r(\,\bar{i}\,) = \overline{z_r(i)}$ and
$$z_r(i) = x(i-1) +1 - i_kp,$$
where $i_k$ is the unique non-negative integer such that $(j-1)p < x(i-1) +1 - i_kp \le jp$ for all $i.$ 
We give an example of $z_r$ below in the case when $p = 3.$ 
\begin{example}
Let $p = 3,$ and let $x = 2.$ 
Then 
$$z_r = (2\ 3)(\bar{2}\ \bar{3})(5\ 6)(\bar{5}\ \bar{6})\ldots(3r-1\ 3r)(\overline{3r-1}\ \overline{3r}),$$
and $(\sigma_1\sigma_2\ldots\sigma_r)^{z_r} = (\sigma_1\sigma_2\ldots\sigma_r)^2.$
\end{example}
For $1 \le i \le r,$ the element $z_r$ commutes with $\tau_i$; moreover, $\sigma_i^{z_r} = \sigma^x_i.$ 
We therefore have the following lemma.
\begin{lemma}\label{lem: norm}
The normaliser subgroup $N_{C_2 \wr S_{rp}}(R_r)$ is generated by the set
$$\{ \tau_i, \sigma_i, \rho_i : 1\le i \le r-1 \} \cup \{\tau_r,\sigma_r\}\cup \{z_r\}.$$
Furthermore, the above set without the element $z_r$ generates the full centraliser $C_{C_2 \wr S_{rp}}(R_r).$
\end{lemma}
It follows that $N_{C_2 \wr S_{rp}}(R_r) \cong (C_{2p} \wr S_r) \rtimes C_{p-1},$ and that $C_{C_2 \wr S_{rp}}(R_r) \cong C_{2p} \wr S_r,$ where the isomorphisms are of abstract groups.
\subsection{The proof of Theorem \ref{thm: main}}
We are now ready to proceed with the first step of the proof. \bigskip \\
\noindent {\bf First step: The Brauer correspondent $M_{(2a,b,c)}(R_r)$.}
Fix $r\in \N$ such that $rp \le n.$
Define
$$T^r = \{(2s,t,u) \in \N_0^3 : 2s+t+u = r, sp \le a, tp \le b, up\le c\}.$$
By the first part of Lemma \ref{lem: ppermbasis}, for each $v \in \mathcal B_{(2a,b,c)}$, there exist $s_v \in \{-1,1\}$ such that
$\{s_v v : v \in \mathcal B_{(2a,b,c)}\}$ is a $p$-permutation basis of $M_{(2a,b,c)}$ with respect to $R_r.$ 
Moreover, by the second part of Lemma \ref{lem: ppermbasis} we can take $s_v =1$ for all $v \in \mathcal B^{R_r}_{(2sp,tp,up)}.$ 
Given $(2s,t,u) \in T^r,$ define $\mathcal A_{(2s,t,u)}$ to be the set
$$\left\lbrace v : \begin{array}{l}
v \in \mathcal B^{R_r}_{(2a,b,c)}\\
\mathcal S(v) \mbox{ contains exactly $2s$ orbits of $\widehat{R_r}$ of length $p$}\\
\mathcal T(v) \mbox{ contains exactly $t$ orbits of $\widehat{R_r}$ of length $p$}\\
\mathcal U(v) \mbox{ contains exactly $u$ orbits of $\widehat{R_r}$ of length $p$}
\end{array}\right\rbrace.$$

\begin{lemma}\label{lem: tensorfact}
There is a direct sum decomposition of $FN_{C_2 \wr S_n}(R_r)$-modules
$$M_{(2a,b,c)}(R_r) = \bigoplus_{(2s,t,u)} \langle \mathcal A_{(2s,t,u)} \rangle,$$
where the sum runs over all $(2s,t,u) \in T^r.$
\end{lemma}
\begin{proof}
Given $v \in \mathcal A_{(2s,t,u)},$ let
$$v = (g, (\{i_{a+1},\overline{i_{a+1}}\},\ldots,\{i_{a+b},\overline{i_{a+b}}\}),([i_{a+b+1},\overline{i_{a+b+1}}],\ldots,[i_n,\overline{i_n}])).$$ 
We first prove that the number of $\widehat{R_r}$-orbits contained in $\mathcal S(v)$ must be even.
If $v \in \mathcal B^{R_r}_{(2a,b,c)},$ then $g \in C_{C_2\wr S_n}(R_r).$ 
Therefore $g$ permutes the $R_r$-orbits as blocks for its action, and the same is true for $\widehat{g}$ and $\widehat{R_r}.$
As $\widehat{g}$ has order 2 and $p$ is odd, the number of $\widehat{R_r}$-orbits contained in $\mathcal S(v)$ is necessarily even. 
Given $h \in N_{C_2 \wr S_n}(R_r),$ let $vh = \pm \tilde v,$ where 
$$\tilde v = (g^h, (\{j_{a+1},\overline{j_{a+1}}\},\ldots,\{j_{a+b},\overline{j_{a+b}}\}),([j_{a+b+1},\overline{j_{a+b+1}}],\ldots,[j_n,\overline{j_n}])).$$ 
The $\widehat{R_r}$-orbits contained in $\mathcal S(\tilde{v})$ are exactly the conjugates, by $\widehat{h}$, of the $\widehat{R_r}$-orbits contained in $\mathcal S(v).$ 
The same argument holds for $\mathcal T(\tilde v)$ and $\mathcal U(\tilde v),$ and so $\tilde v \in \langle\mathcal  A_{(2s,t,u)} \rangle.$ 
It follows that $\langle\mathcal A_{(2s,t,u)} \rangle$ is a submodule of $M_{(2a,b,c)}(R_r).$ 
The lemma now follows as $\mathcal B_{(2a,b,c)} =  \bigcup \mathcal A_{(2s,t,u)}.$
\end{proof}
In the following lemma, we factorise the module $\langle \mathcal A_{(2s,t,u)}\rangle$ as an outer tensor product of modules, compatible with the factorisation of $N_{C_2 \wr S_n}(R_r)$ in (\ref{eqn: normaliser_fact}).
By doing this, we see that in order to understand $M_{(2a,b,c)}(R_r),$ it is sufficient to understand the modules $M_{(2sp,tp,up)}(R_r),$ where $(2s,t,u) \in T^r.$ 

\begin{lemma}\label{lem: tensorfact1}
There is an isomorphism
$$\langle \mathcal A_{(2s,t,u)} \rangle \cong M_{(2sp,tp,up)}(R_r) \boxtimes M_{(2(a-sp),b-tp,c-up)},$$
of $F[N_{C_2 \wr S_{rp}}(R_r) \times C_2\wr S_{\{rp+1,\ldots,n\}}]$-modules. 
\end{lemma} 
\begin{proof}
Let $\mathcal B_{(2(a-sp),b-tp,c-up)}^+$ be the set of elements in $\mathcal B_{(2(a-sp),b-tp,c-up)},$ each shifted by appropriately by $rp$ or $\overline{rp}.$
The $F$-span of $\mathcal B^+_{(2(a-sp),b-tp,c-up)}$ is therefore an $F[C_2 \wr S_{\{rp+1,\ldots,n\}}]$-module isomorphic to $M_{(2(a-sp),b-tp,c-up)}.$

Let $v \in \mathcal A_{(2s,t,u)}$ be such that
$$v = (g, (\{i_{a+1},\overline{i_{a+1}}\},\ldots,\{i_{a+b},\overline{i_{a+b}}\}),([i_{a+b+1},\overline{i_{a+b+1}}],\ldots,[i_n,\overline{i_n}])),$$
where $\mathcal S(v) = \{i_1,\ldots,i_a\}$ and the notation is chosen so that
$$\{i_1, \ldots, i_{2sp}\} \cup \{i_{a+1},\ldots,i_{a+tp}\} \cup \{i_{a+b+1},\ldots,i_{a+b+up}\}  = \{1,2,\ldots,rp\}.$$
Let $v_1 \in \mathcal{B}^{R_r}_{(2sp,tp,up)}$ be the unique element such that 
\begin{align*}
\mathcal S(v_1) &= \mathcal S(v) \cap \{1,2,\ldots,rp\}\\
\mathcal T(v_1) &= \mathcal T(v) \cap \{1,2,\ldots,rp\}\\
\mathcal U(v_1) &= \mathcal U(v) \cap \{1,2,\ldots,rp\}.
\end{align*}
By construction, the $p$-element $\sigma_1\sigma_2\ldots\sigma_{r}$ has support $\{1,\overline{1},\ldots,rp,\overline{rp}\},$ and so $v$ is fixed by $R_r$ if and only if $v_1$ is fixed by $R_r$. 
Let $v_2 \in \mathcal B^+_{(2(a-sp),b-tp,c-up)}$ be such that
\begin{align*}
\mathcal S(v_2) &= \mathcal S(v)\backslash \mathcal S(v_1) \\
\mathcal T(v_2) &=  \mathcal T(v)\backslash \mathcal T(v_1) \\
\mathcal U(v_2) &=  \mathcal U(v)\backslash \mathcal U(v_1) .
\end{align*}
It follows that there is a natural bijection $f$ between $\mathcal{B}^{R_r}_{(2a,b,c)}$ and 
$$\mathcal{B}^{R_r}_{(2sp,tp,up)} \times \mathcal B_{(2(a-sp),b-tp,c-up)}^+,$$
defined by  $f(v) = v_1 \otimes v_2.$ 

We now show that $f$ is an $F[N_{C_2\wr S_{n}}(R_r)]$-module homomorphism.
Given $g \in N_{C_2\wr S_{n}}(R_r),$ and $v \in \mathcal B_{(2a,b,c)},$ let $v^\star \in \mathcal B_{(2a,b,c)}$ be such that the entries in its $b$-tuple and $c$-tuple are those of $vg$ in ascending order (with respect to the orders in \textsection \ref{sec: hypspecht}).
Let $h \in C_2 \wr S_{\mathcal T(vg)} \times C_2 \wr S_{\mathcal U(vg)}$ be the unique permutation such that $v^\star = vgh.$
As $g \in N_{C_2\wr S_{rp}}(R_r) \times C_2 \wr S_{\{rp+1,\ldots,n\}},$ it permutes the elements in the sets $\{1,2,\ldots,rp,\bar{1},\ldots,\overline{rp}\}$ and $\{rp+1,\ldots,n,\overline{rp+1},\ldots,\bar{n}\}$ separately. 
It follows that there exists a factorisation $h = h_1h_2,$ where $h_1 \in C_2 \wr S_{\{1,2,\ldots,rp\}}$ and $h_2 \in C_2 \wr S_{\{rp+1,\ldots,n\}}.$
Therefore
\begin{align*}
f(vg) &= f(\sgn(\,\widehat h\,) v^\star )\\
&= \sgn(\,\widehat h\,)(v^\star_1 \otimes v^\star_2)\\
&= \sgn(\,\widehat h\,)\sgn(\,\widehat{h_1}\,)\sgn(\,\widehat{h_2}\,)(v_1\otimes v_2)\\
&= f(v)g,
\end{align*}
and so the result is proved.
\end{proof}

In order to express $M_{(2sp,tp,up)}(R_r)$ as a sum of indecomposable modules,
we first write $M_{(2sp,tp,up)}(R_r)$ as a direct sum of $FN_{C_2\wr S_{rp}}(R_r)$-modules $N_{(\lambda,t,u)}$ (defined below), before showing that each of these modules is indecomposable. 
We require a deeper understanding of the fixed points $v \in \mathcal B^{R_r}_{(2sp,tp,up)}$ before we can define $N_{(\lambda,t,u)}.$
We begin by considering the example $M_{(2p,0,0)}(R_r)$ for all $r \in \N.$ 
This illustrative example will be used when describing $\mathcal B^{R_r}_{(2sp,tp,up)}$ in the general case. 
\begin{example}\label{ex: case_1}
The $FC_2 \wr S_{2p}$-module $M_{(2p,0,0)}$ is a permutation module, with permutation basis given by the set
$$\mathcal B_p := \{f_p^h : h \in C_2 \wr S_{2p}\}.$$
The set $T^j$ is empty for all $j\in \N$ such that $j\neq 2,$ and so we consider $M_{(2p,0,0)}(R_2).$ 
We write $\sigma_1\sigma_2$ as
$$\sigma_1\sigma_2 = (1\ 2\ldots\ p)(\,\overline{1}\ \overline{2}\dots\ \overline{p}\,)(1^*\ 2^*\ldots\ p^*)(\,\overline{1^*}\ \overline{2^*}\dots\ \overline{p^*}\,),$$
where $x^* := x+p$ for $1 \le x \le p.$  

Let $g \in\mathcal B_p$ be fixed by $R_2.$ 
If $(1)g = x,$ then $(2)g = (x)\sigma_1\sigma_2.$
Therefore, for $2 \le i \le p-1,$
\begin{equation}\label{eqn: case_1} 
(i+1)g = (x)(\sigma_1\sigma_2)^i,
\end{equation}
and so $g$ is completely determined by $(1)g.$

Suppose that $x \in \{2,\overline{2},\ldots,p,\overline{p}\}.$ 
If $x\in \{2,\ldots,p\},$ it follows from (\ref{eqn: case_1}) that $(x)g = 2x -1 \mod p.$ 
As $p$ is odd, we cannot have that $(x)g = 1,$ and so $g$ does not have order 2. 
It follows that $g$ cannot be a conjugate of $f_p,$ which is a contradiction. 
A similar argument shows that $x \not\in \{\overline{2},\ldots,\overline{p}\}.$ 

There are now $2p$ possible choices for $x.$ 
As each such choice completely determines the permutation $g,$ the module $M_{(2p,0,0)}(R_2)$ has dimension $2p.$ 
\end{example}

Fix $(2s,t,u) \in T^r,$ and let $k = t+u.$
We define $\Omega^{(2s;k)}$ to be the set of elements of the form
$$\{\{i_1,i'_1\},\ldots,\{i_s,i'_s\},\{j_1,\ldots,j_k\}\},$$
where $\{i_1,i'_1,\ldots,i_s,i'_s,j_1,\ldots,j_k\} = \{1,2,\ldots,r\}.$
Let $c_{s,k} = \vert \Omega^{(2s;k)} \vert.$ 
Given $\omega \in \Omega^{(2s;k)}$ of the above form, define 
$$R_\omega = \langle\sigma_{i_1}\sigma_{i'_1}\rangle \times \cdots \times \langle \sigma_{i_s}\sigma_{i'_s}\rangle \times \langle \sigma_{j_1}\rangle \times \cdots \times \langle \sigma_{j_k} \rangle.$$
and write $\mathcal B(\omega)$ for $\mathcal B^{R_\omega}_{(2sp,tp,up)}.$
\begin{lemma}\label{lem: omega_type}
Let $v \in \mathcal B^{R_r}_{(2sp,tp,up)}.$
Then $v \in \mathcal B(\omega)$ for a unique $\omega \in \Omega^{(2s;k)}.$ 
\end{lemma}
\begin{proof}
By the second part of Lemma \ref{lem: ppermbasis}, the vector $v \in \mathcal B_{(2sp,tp,up)}$ is fixed by $R_r$ if and only if $\theta(v)$ is fixed by $R_r.$
Let $v$ be such that $\theta(v) = (g,\gamma,\delta)$ where
\begin{align*}
\gamma&=\{\{i_{2sp+1},\overline{i_{2sp+1}}\},\ldots,\{i_{(2s+t)p},\overline{i_{(2s+t)p}}\}\}\\
\delta &= \{[i_{(2s+t)p+1},\overline{i_{(2s+t)p+1}}],\ldots,[i_{rp},\overline{i_{rp}}]\}.
\end{align*}
By definition, $g$ can be written a product of $s$ disjoint $p$ positive $2$-cycles $g_1,\ldots,g_s.$
Let $\{i_1,i_2,\ldots,i_{2s-1},i_{2s}\}$ be such that $\supp(\sigma_{i_{2j-1}}\sigma_{i_{2j}}) = \supp(g_j)$ for each $j.$
It follows from Example \ref{ex: case_1} that $g$ commutes with $R_r$ if and only if $g$ commutes with $\sigma_{i_{2j-1}}\sigma_{i_{2j}}$ for each $j.$

Let $\{j_1,\ldots,j_t\}$ be such that 
$$\supp(\sigma_{j_1}\ldots\sigma_{j_t}) = \{i_{2sp+1},\overline{i_{2sp+1}},\ldots,i_{(2s+t)p},\overline{i_{(2s+t)p}}\}.$$
As $\gamma$ is fixed by $R_r,$ the set $\mathcal T(v)$
is equal to a union of $\widehat{R_r}$-orbits. 
The orbits of $\widehat{R_r}$ are equal to precisely the orbits of $\widehat{\sigma_{j_i}}$ for each $1 \le i \le r.$ 
Therefore $\gamma$ is fixed by $R_r$ if and only if it is fixed by the group
$\langle \sigma_{j_1}\rangle \times \cdots \times \langle \sigma_{j_t}\rangle.$ 

Similarly if $\delta$ is such that 
$$\supp(\sigma_{k_1}\ldots\sigma_{k_u}) = \{i_{(2s+t)p+1},\overline{i_{(2s+t)p+1}},\ldots,i_{rp},\overline{i_{rp}}\},$$ 
then $\delta$ is fixed by the group $\langle \sigma_{k_1}\rangle \times \cdots \times \langle \sigma_{k_u}\rangle.$

Therefore if $v$ is fixed by $R_r,$ then $v$ is fixed by $R_\omega$, where 
$$\omega = \{\{{i_1},{i_2}\},\ldots,\{{i_{2s-1}},{i_{2s}}\},\{j_1,\ldots,j_t,k_1,\ldots,k_u\}\}.$$
Moreover, $\omega$ is unique as it is determined by the fixed sets $\supp(g),\supp(\gamma),$ and $\supp(\delta).$
\end{proof}

Given $\varnothing \neq E \subseteq \{1,2,\ldots,r\},$ we define $\tau_E = \prod_{e\in E}\tau_e.$ If $E$ is empty, then we define $\tau_E = 1.$
\begin{definition}
Fix $y \in \{-1,1\}^s.$
Given $(g,\gamma,\delta) \in \mathcal B_{(2sp,tp,up)}^{R_r},$ let 
$$\omega :=\{\{i_1,i'_1\},\ldots,\{i_s,i'_s\},\{j_1,\ldots,j_k\}\}$$
be the unique element of $\Omega^{(2s;k)}$ such that $(g,\gamma,\delta) \in \mathcal B(\omega).$
Define
$$(y(g),\gamma,\delta) = \sum_{E\subseteq\{i_1,\ldots,i_s\}}(\prod_{e\in E}y_e)(g^{\tau_E},\gamma,\delta).$$
It follows from Example \ref{ex: case_1} and Lemma \ref{lem: omega_type} that $\tau_{i_j}$ and $\tau_{i'_j}$ act in the same way on $g,$ and so $(y(g),\gamma,\delta)$ is well-defined. 
\end{definition}
Given $y\in \{-1,1\}^s$, let $\lambda \in \Lambda(2,s)$ be the composition of $s$ such that $\lambda_1$ (resp. $\lambda_2$) is equal to the number of $y_i$ that are equal to $+1$ (resp. $-1$). 
We say that $y$ has {\it weight} $\lambda$ in this case. 

We now define $N_{(\lambda,t,u)}$ to be the $F$-span of
\begin{equation}\label{eq: basisNlambda}
\{(y(g),\gamma,\delta): (g,\gamma,\delta)\in \mathcal B^{R_r}_{(2sp,tp,up)} \mbox{ and $y$ has weight } \lambda\}.
\end{equation}
It is clear that 
$$M_{(2sp,tp,up)}(R_r) = \bigoplus_{\lambda \in \Lambda(2,s)} N_{(\lambda,t,u)},$$
is an equality of vector spaces. 

Before we state the next result, we remark that $N_{C_2 \wr S_{rp}}(R_r)$ permutes the $R_r$-orbits as blocks for its action. 
It follows that the set of subgroups of the form $R_\omega$ is an $N_{C_2 \wr S_{rp}}(R_r)$-set, with the action given by conjugation.  
We have seen in the proof of Lemma \ref{lem: omega_type} that the $R_r$-orbits are the same as the orbits of the subgroup $C := \langle \sigma_1 \rangle \times \cdots \times \langle \sigma_r \rangle,$ and we write $\mathcal O_i$ for the union of the non-trivial orbits of $\langle \sigma_i \rangle.$ 
\begin{lemma}\label{lem: tau_action}
Given $\omega,\widetilde{\omega} \in\Omega^{(2s;k)},$ let $h \in N_{C_2\wr S_{rp}}(R_r),$ be such that $R_\omega^h = R_{\widetilde{\omega}}.$
Given $1 \le i \le r,$ let $\tilde{i}$ be such that $\mathcal O_i^h = \mathcal O_{\tilde{i}}.$ 
Then $(y(g),\gamma,\delta)h$ is contained in the $F$-span of $\mathcal B(\,\widetilde\omega\,).$ 
\end{lemma}
\begin{proof}
It follows from the definition of $(y(g),\gamma,\delta)$ that
\begin{align*}
(y(g),\gamma,\delta){h} &= \sum_{E\subseteq\{i_1,\ldots,i_s\}}(\prod_{e\in E}y_e)(g^{\tau_E},\gamma,\delta)h\\ 
&= \sum_{E\subseteq\{i_1,\ldots,i_s\}}(\prod_{e\in E}y_e)((g^{\tau_E})^h,\gamma{h},\delta{h})\\
&= \sum_{\widetilde{E}\subseteq\{\widetilde{i_1},\ldots,\widetilde{i_s}\}}(\prod_{e\in E}y_e)((g^h)^{\tau_{\widetilde{E}}},\gamma{h},\delta{h})\\
&= (\widetilde{y}(g^h),\gamma h, \delta h),
\end{align*}
where $\widetilde{E} = \{\tilde{i} : i \in E\}$ and $\widetilde y_{\tilde{i}} = y_i$ for all $i \in \{i_1,i_2,\ldots,i_s\}.$
The lemma is proved once we show that $({g}^h,\gamma{h},\delta{h})$ is fixed by $R_{\widetilde{\omega}}.$ 
As $\sigma_i^h = \sigma_{\tilde{i}}$ for all $1 \le i \le r,$ for $1\le j \le s$
$$(g^{h})^{\sigma_{\widetilde{i_j}}\sigma_{\widetilde{i'_j}}} = (g^{(\sigma_{i_j}\sigma_{i'_j})})^{h} = g^h.$$
An entirely similar argument shows that $\gamma{h}\sigma_{\widetilde{i_j}} = \gamma{h}$ for $s < j \le s+t$, and that $\delta{h}\sigma_{\widetilde{i_j}} = \delta{h}$ for $s+t < j \le r.$ 
\end{proof}
It follows from Lemma \ref{lem: tau_action} that each $N_{(\lambda,t,u)}$ is an $FN_{C_2\wr S_{rp}}(R_r)$-module. 
\begin{corollary}\label{rem: tau_special}
Let $h \in N_{C_2\wr S_{rp}}(R_r)$ be such that $\tau_i^h = \tau_i$ for $1 \le i \le r,$ and $\sigma_i^h = \sigma_i^x,$ for some $x\in \N.$
If $(g,\gamma,\delta) \in \mathcal B(\omega),$ then $(g,\gamma,\delta)h$ is contained in the $F$-span of $B(\omega).$ 
In particular if $h = \tau_{i_j}$ for some $j \in \{1,2,\ldots,s\},$ then $(y(g),\gamma,\delta)h = y_j(y(g),\gamma,\delta).$
\end{corollary}
\begin{proof}
For the first statement, observe that $\mathcal O_i^h = \mathcal O_i$ for $1 \le i \le r.$  
Now apply Lemma \ref{lem: tau_action}.  
For the second statement observe that when $h = \tau_{i_j}$ for some $j \in \{1,\ldots,s\},$ then
\begin{align*}
(y(g),\gamma,\delta){\tau_{i_j}} &= \sum_{E\subseteq\{i_1,\ldots,i_s\}}(\prod_{e\in E}y_e)(g^{\tau_E\tau_{i_j}},\gamma,\delta) \\
&= y_j\sum_{E\subseteq\{i_1,\ldots,i_s\}}(\prod_{e\in E}y_ey_j)(g^{\tau_E\tau_{i_j}},\gamma,\delta) = y_j(y(g),\gamma,\delta).
\end{align*}
\end{proof}

Write $K_r$ for $C_{C_2 \wr S_{rp}}(R_r),$ which recall is isomorphic to $C_{2p} \wr S_r.$ 
In order to prove that each $N_{(\lambda,t,u)}$ is an indecomposable $FN_{C_2 \wr S_{rp}}(R_r)$-module, we show that it is indecomposable as an $FK_r$-module. 
We do this by filling in the details of the following sketch.

Given $1\le i \le r,$ define $D_i = \langle \sigma_i, \tau_i \rangle,$ and so $D_1 \times \cdots \times D_r$ is a normal subgroup in $K_r.$
We define an $F[D_1 \times \cdots \times D_r]$-module $N^{\omega^\star}_y,$ and in Lemma \ref{lem: inertial} we determine its inertial group $Y_{(\lambda,t,u)}$ in $K_r.$ 
Using Lemma \ref{lem: case_r}, we determine the dimension of $N_{(\lambda,t,u)}.$
In Lemma \ref{lem: R2r_decomp}, we first show that $N_y:=N^{\omega^\star}_y\Ind_{X_{(\lambda,t,u)}}^{Y_{(\lambda,t,u)}}$ is indecomposable, where $X_{(\lambda,t,u)}$ is the largest subgroup in $Y_{(\lambda,t,u)}$ that $N^{\omega^\star}_y$ can be extended to. 
We then prove that $N_{(\lambda,t,u)} = N_y\Ind_{Y_{(\lambda,t,u)}}^{K_r}.$
It follows using Clifford theory (see \cite[Proposition 3.13.2]{Benson}) that $N_{(\lambda,t,u)}$ is an indecomposable $FK_r$-module. 

Define $\omega^{\star} = \{\{1,s+1\},\ldots,\{s,2s\},\{2s+1,\ldots,r\}\} \in \Omega^{(2s;k)}.$
Furthermore define
$v^\star=(f_{sp},\gamma^\star,\delta^\star) \in B(\omega^\star),$
where
\begin{align*}
\mathcal T(v^\star) &= \supp(\sigma_{2s+1}\ldots\sigma_{2s+t})\cap \{1,2,\ldots,n\}\\
\mathcal U(v^\star) &= \supp(\sigma_{2s+t+1}\ldots\sigma_r)\cap \{1,2,\ldots,n\}.
\end{align*}
Given $\lambda \in \Lambda(2,s),$ define $y_\lambda \in \{-1,1\}^s$ to be the tuple of weight $\lambda$ such that 
$$(y_\lambda)_i = \left\lbrace \begin{array}{rl}
1 & \mbox{if $1 \le i \le \lambda_1$}\\
-1 &\mbox{if $\lambda_1 + 1 \le i \le s$}.
\end{array}\right.$$
Define $N^{\omega^\star}_y$ to be the $F$-span of 
$$\{(y_\lambda(g),\gamma^\star,\delta^\star) : (g,\gamma^\star,\delta^\star)\in B(\omega^\star)\}.$$
We also define $X_{(\lambda,t,u)}$ to be the subgroup of $K_r$ generated by the set
\begin{align*}
\{\sigma_i,\tau_i : 1 \le i \le r\} \cup \{\rho_1^{\rho_2\rho_3\ldots\rho_s}\} &\cup \{\rho_i\rho_{i+s} : 1 \le i \le s-1 \mbox{ and } i \neq \lambda_1\}\\ 
&\cup \{\rho_i : 2s+1 \le i < r, i \neq 2s+t\} ,
\end{align*}
and $Y_{(\lambda,t,u)}$ to be the subgroup of $K_r$ generated by the set
$$\{\sigma_i,\tau_i : 1 \le i \le r\} \cup \{\rho_i : 1\le i \le r-1 \mbox{ and } i \not\in \{2\lambda_1, 2s,2s+t\}\}.$$
Similar to the remark following Lemma \ref{lem: norm}, there are isomorphisms of abstract groups $X_{(\lambda,t,u)} \cong C_{2p}\wr ((S_2\wr S_{\lambda}) \times S_{t}\times S_{u}),$ and $Y_{(\lambda,t,u)} \cong C_{2p} \wr (S_{2\lambda} \times S_t \times S_u).$
\begin{lemma}\label{lem: inertial}
The vector space $N^{\omega^\star}_y$ is an $F[D_1 \times \cdots \times D_r]$-module, with inertial group $Y_{(\lambda,t,u)}$ in $K_r.$
Moreover, we can extend $N^{\omega^\star}_y$ to a module for $FX_{(\lambda,t,u)}.$
\end{lemma}
\begin{proof}
That $N^{\omega^\star}_y$ is an $F[D_1 \times \cdots \times D_r]$-module follows from by applying the first statement of Corollary \ref{rem: tau_special}. 

Write $T$ for the inertial group of $N^{\omega^\star}_y,$ which permutes the groups $D_i$ by conjugation.
The permutations $\sigma_1, \ldots, \sigma_{2s}$ act freely on $N^{\omega^\star}_{y},$ whereas $\sigma_{2s+1},\ldots,\sigma_r$ all act trivially on $N^{\omega^\star}_y.$ 
Therefore $T$ must be contained in the subgroup of $K_r$ that permutes the groups $D_1,\ldots,D_{2s}$ amongst themselves and the groups $D_{2s+1},\ldots,D_{r}$ amongst themselves. 

For $2s < i \le 2s+t,$ the action of $\tau_i$ on $(y_\lambda(g),\gamma^\star,\delta^\star)$ is determined by its action on 
$(\{(i-1)p+1,\overline{(i-1)p+1}\},\ldots,\{ip,\overline{ip}\}).$ 
Therefore $\tau_i$ acts trivially in this case.
Similarly for $2s+t < i \le r,$ the action of $\tau_i$ on $(y_\lambda(g),\gamma^\star,\delta^\star)$ is determined by its action on $([(i-1)p+1,\overline{(i-1)p+1}\,],\ldots,[ip,\overline{ip}]).$
It follows that $\tau_i$ acts with sign $(-1)^p,$ which is negative as $p$ is odd. 
Therefore $T$ must be contained in the subgroup of $K_r$ that permutes the $D_{2s+1},\ldots,D_{2s+t}$ amongst themselves, and the $D_{2s+t+1},\ldots,D_r$ amongst themselves.

It follows from the second statement of Corollary \ref{rem: tau_special} that $T$ must permute the groups $D_1,\ldots,D_{\lambda_1},D_{s+1},\ldots,D_{\lambda_1+s}$ amongst themselves, and the same is true for the groups $D_{\lambda_1+1},\ldots,D_{2s},D_{\lambda_1+s+1},\ldots,D_{2s}.$
This shows that $T$ must be contained in $Y_{(\lambda,t,u)}.$
Moreover, if $h \in Y_{(\lambda,t,u)},$ then $(N^{\omega^\star}_y)^h \cong N^{\omega^\star}_y.$ 
Therefore $Y_{(\lambda,t,u)}$ is contained in $T,$ which proves the second statement of the lemma. 

For the final statement, it remains to prove that $N^{\omega^\star}_y$ is closed under the action of 
$$Z \cup \{\rho_i : 2s+1 \le i < r, i \neq 2s+t\},$$
where $Z = \{\rho_1^{\rho_2\rho_3\ldots\rho_s}\} \cup \{\rho_i\rho_{i+s} : 1 \le i \le s-1 \mbox{ and } i \neq \lambda_1\}.$ 
It is sufficient to prove that each of $(y(g)^{z},\gamma^\star,\delta^\star)$, where $z \in Z$, and  $(y(g),(\gamma^\star){\rho_i},\delta^\star),$ where $2s+1 \le i < 2s+t$, and $(y(g),\gamma^\star, (\delta^\star)\rho_i),$ where $2s+t < i < r,$ is contained in $N^{\omega^\star}_y.$ 

First consider $\gamma^\star{\rho_i},$ where $2s+1 \le i < 2s+t.$
As $\rho_i$ permutes precisely those orbits of $R_{\omega^\star}$ with support equal to the support of $\gamma^\star,$ it follows that
$\gamma^\star{\rho_i} = \pm \gamma^\star.$
The same argument shows that $\delta^\star{\rho_i} = \pm \delta^\star$ for $2s+t < i < r.$
If $z \in \{\rho_1^{\rho_2\rho_3\ldots\rho_s}\} \cup \{\rho_i\rho_{i+s} : 1 \le i \le s-1 \mbox{ and } i \neq \lambda_1\},$ then
\begin{align*}
(y(g),\gamma^\star,\delta^\star){z} &= \sum_{E\subseteq\{1,\ldots,s\}}(\prod_{e\in E}y_e)(g^{(\tau_E)z},\gamma^\star,\delta^\star)\\
&= \sum_{E'\subseteq\{1,\ldots,s\}}(\prod_{e\in E}y_e)((g^{z})^{\tau_{E'}},\gamma^\star,\delta^\star),
\end{align*}
where $E' = E$ if $z =\rho_1^{\rho_2\rho_3\ldots\rho_s}$, otherwise $E'$  is the subset of $\{1,\ldots,s\}$ obtained from $E$ by swapping $i$ and $i+1.$ 
As $i \neq \lambda_1,$ we have $y_{i} = y_{i+1}$ in all cases. 
It follows that $(y(g),\gamma,\delta){z} = (y(g^{z}),\gamma,\delta).$
The lemma is proved if $(g^{z},\gamma,\delta) \in B(\omega^\star).$ 
This follows from the first statement of Corollary \ref{rem: tau_special} as $z$ centralises $R_{\omega^\star}.$
\end{proof}

\begin{lemma}\label{lem: case_r}
The module $M_{(2sp,tp,up)}(R_r)$ has dimension equal to
$$(2p)^s \times \binom{k}{t} \times c_{s,k}.$$  
\end{lemma}
\begin{proof}
By Lemma \ref{lem: omega_type} every element in $\mathcal B_{(2sp,tp,up)}^{R_r}$ is fixed by $R_\omega,$ for a unique $\omega \in \Omega^{(2s;k)}.$
We therefore count the size of $\mathcal B(\omega)$ for each $\omega.$ 
Fix $\omega \in\Omega^{(2s;k)},$ and let $\omega =\{\{i_1,i'_1\},\ldots,\{i_s,i'_s\},\{j_1,\ldots,j_k\}\}$. 

Let $(g,\gamma,\delta) \in \mathcal B(\omega)$ be such that $g = g_1\ldots g_s,$ where $g_j$ is fixed by $\sigma_{i_j}\sigma_{i'_j}$ for each $j.$
By Example \ref{ex: case_1}, each $\sigma_{i_j}\sigma_{i'_j}$ has $2p$ fixed points in $\mathcal B_{(2p,0,0)}$. 
Therefore there are $(2p)^s$ choices for the fixed points of
$\langle \sigma_{i_1}\sigma_{i'_1}\rangle \times \cdots \times \langle \sigma_{i_s}\sigma_{i'_s}\rangle$
under the conjugacy action.

Let $\gamma :=(\{\gamma_1,\overline{\gamma_1}\},\{\gamma_2,\overline{\gamma_2}\},\ldots,\{\gamma_{tp},\overline{\gamma_{tp}}\})$
be such that $\gamma_1< \gamma_2 < \cdots <\gamma_{tp}$ and
$\supp(\sigma_{j_1}\ldots\sigma_{j_t}) = \{\gamma_1,\overline{\gamma_1},\ldots,\gamma_{tp},\overline{\gamma_{tp}}\}.$
Then $\gamma$ is the unique element of this form with support not disjoint to $\sigma_{j_1}\ldots\sigma_{j_t}$ that is fixed by $\sigma_{j_1}\ldots\sigma_{j_t}.$ 
Similarly, we define $\delta = ([\delta_1,\overline{\delta_1}],\ldots,[\delta_{up},\overline{\delta_{up}}])$ to be such that $\delta_1 < \delta_2 < \cdots < \delta_{up}$ and $\{\delta_1,\overline{\delta_1},\ldots,\delta_{tp},\overline{\delta_{tp}}\} = \supp(\sigma_{j_{t+1}}\ldots\sigma_{j_k}).$
Then $\delta$ is the unique element with support not disjoint to $\sigma_{j_{t+1}}\ldots\sigma_{j_k}$ that is fixed by $\sigma_{j_{t+1}}\ldots\sigma_{j_k}.$

As there are $\binom{k}{t}$ ways to choose $j_1,j_2,\ldots,j_t,$ there are $(2p)^s \times \binom{k}{t}$ fixed points of $R_\omega$ in $\mathcal B_{(2sp,tp,up)}.$
The statement of the lemma now follows by definition of $c_{s,k}.$
\end{proof}
Recall that $C$ is defined to be the elementary abelian group $\langle \sigma_1,\ldots, \sigma_r\rangle.$
\begin{lemma}\label{lem: R2r_decomp}
The module $N_{(\lambda,t,u)}$ is an indecomposable $FK_r$-module. 
\end{lemma}
\begin{proof}
Define $\Omega^{(2\lambda;k)}$ to be the subset of $\Omega^{(2s;k)}$ consisting precisely of the $\omega \in \Omega^{(2s;k)}$ of the form
$\{\{i_1,i_1'\},\ldots,\{i_s,i_s'\},\{j_1,\ldots,j_k\}\}$ such that 

\begin{align*}
\{i_1,i'_1,\ldots,i_{\lambda_1},i'_{\lambda_1}\} &= \{1,\ldots,\lambda_1,s+1,\ldots,\lambda_1+s\}\\
\{i_{\lambda_1+1},i'_{\lambda_1+1},\ldots,i_{s},i'_{s}\} &= \{\lambda_1+1,\ldots,s,\lambda_1+1+s,\ldots,2s\}\\
\{j_1,\ldots,j_k\} &= \{2s+1,\ldots,r\},
\end{align*}
and let $c_{\lambda,k} = \size{\Omega^{(2\lambda;k)}}.$ 
The module $N_{y}:=N^{\omega^\star}_{y}\big\uparrow_{X_{(\lambda,t,u)}}^{Y_{(\lambda,t,u)}}$ has a basis given by the set 
$$\{(y(g),\gamma^\star,\delta^\star) : (g,\gamma^\star,\delta^\star) \in \mathcal B(\omega), \omega \in \Omega^{(2\lambda;k)}\}.$$
Therefore $N_{y}(R_{\omega^\star})$ and $N^{\omega^\star}_{y}$ are equal as vector spaces.
By the second paragraph in the proof of Lemma \ref{lem: case_r}, there are $(2p)^s$ choices for $g$ in $(g,\gamma^\star,\delta^\star).$ 
Given $E \subseteq \{i_1,\ldots,i_s\},$ the second statement of Corollary \ref{rem: tau_special} implies that $(y(g),\gamma^\star,\delta^\star)$ and $(y(g),\gamma^\star,\delta^\star)\tau_E$ are equal up to a sign. 
There are $2^s$ choices for $E,$ and so $N^{\omega^\star}_y$ has dimension $p^s.$

The group $R_{\omega^\star}$ acts trivially on $N^{\omega^\star}_y,$ and so by \cite[Chapter 8, Corollary 3]{Alperin}
$$N_{y}(R_{\omega^\star})\Res_C \cong F\Ind_{R_{\omega^\star}}^C,$$
where the subgroup $C$ is defined just before this lemma.
Then $N_y(R_{\omega^\star})$ is an indecomposable $FC$-module, and so $N_{y}(R_{\omega^\star})$ is an indecomposable $F[N_{Y_{(\lambda,t,u)}}(R_{\omega^\star})]$-module.
It follows that there exists a unique summand of $N_{y}$ with vertex containing $R_{\omega^\star}.$ 
Let $W$ be a non-zero indecomposable summand of $N_y.$ 
As
$$N_{y}\Res_C \cong \bigoplus_{\omega \in \Omega^{(2\lambda;k)}} F\Ind_{R_{\omega}}^C,$$
the Krull--Schmidt theorem implies that each indecomposable summand of $W\hspace{-0.5pt}\res_C$ is isomorphic to $F\ind_{R_{\omega^\star}}^C.$
Therefore $W(R_{\omega^\star}) \neq 0,$ and so Lemma \ref{lem: brauer_vertex} states that $W$ has a vertex containing $R_{\omega^\star}.$ 
The module $N_{y}$ is therefore indecomposable. 

Let $(\widetilde{y}(g),\tilde{\gamma},\tilde{\delta}) \in N_{(\lambda,t,u)}$ be such that $(g,\tilde{\gamma},\tilde{\delta}) \in B(\widetilde{\omega}).$ 
As $\widetilde{y}$ has weight $\lambda,$ and $K_r$ permutes the $R_r$-orbits transitively, it follows from Lemma \ref{lem: tau_action} that there exists $\rho \in \langle \rho_1,\ldots,\rho_{r-1} \rangle$
such that $\pm (\widetilde{y}(g),\tilde{\gamma},\tilde{\delta}) = (y(g^{\rho^{-1}}),\gamma,\delta)\rho,$
where $(y(g^{\rho^{-1}}),\gamma,\delta) \in N_y.$ 
Therefore $N_y$ generates $N_{(\lambda,t,u)}$ as an $F{K_r}$-module.

By definition there are $c_{s,k}$ choices for $\omega \in \Omega^{(2s;k)},$ and there are $\binom{s}{\lambda_1}$ choices for $y \in \{-1,1\}^s$ of weight $\lambda.$ 
Therefore $N_{(\lambda,t,u)}$ has dimension $c_{s,k} \times \binom{s}{\lambda_1} \times p^s \times \binom{k}{t}.$ 
As $N_y$ has dimension $c_{\lambda,k} \times p^s,$ applying \cite[Chapter 8, Corollary 3]{Alperin} gives
$$N_{(\lambda,t,u)}\cong N_{y}\big\uparrow_{Y_{(\lambda,t,u)}}^{K_r}.$$
Lemma \ref{lem: inertial} states that $Y_{(\lambda,t,u)}$ is the inertial group of the $F[D_1\times \cdots \times D_r]$-module $N^{\omega^\star}_y.$ 
As $N^{\omega^\star}_y$ is extended from $D_1\times \cdots \times D_r$ to $X_{(\lambda,t,u)},$ we have that
$N_y\Res_{D_1\times \cdots \times D_r}$ is isomorphic to a direct sum of $[Y_{(\lambda,t,u)} : X_{(\lambda,t,u)}]$ copies of $N^{\omega^\star}_y.$
Therefore the proof of Proposition 3.13.2 in \cite{Benson} carries over to this case, and so $N_{(\lambda,t,u)}$ is an indecomposable $FK_r$-module. 
\end{proof}
By Lemma \ref{lem: broue_norm}, $M_{(2sp,tp,up)}(R_r)$ and $N_{(\lambda,t,u)}$ are $p$-permutation $FN_{C_2 \wr S_{rp}}(R_r)$-modules, where $\lambda \in \Lambda(2,s).$ 
Write $J_r$ for $N_{C_2 \wr S_{rp}}(R_r).$ 
As $R_r \unlhd R_{{\omega^\star}},$ we have $R_{\omega^\star} \le J_r.$ 
By Lemma \ref{lem: broue_norm}
$M_{(2sp,tp,up)}(R_r)(R_{\omega^\star}) \cong M_{(2sp,tp,up)}(R_{\omega^\star}),$
where the isomorphism is of $FN_{J_r}(R_{\omega^\star})$-modules. 
Then Lemma \ref{lem: R2r_decomp} implies that 
$$M_{(2sp,tp,up)}(R_{\omega^\star}) \cong \bigoplus_{\lambda \in \Lambda(2,s)} N_{(\lambda,t,u)}(R_{\omega^\star}),$$
as $FN_{J_r}(R_{\omega^\star})$-modules.
Moreover, for all $\lambda \in \Lambda(2,s),$ the basis defining $N_{(\lambda,t,u)}$ in \eqref{eq: basisNlambda} is a $p$-permutation basis of $N_{(\lambda,t,u)}$ with respect to $R_{\omega^\star}.$ 

Recall that $U$ is a non-projective indecomposable summand of $M_{(2a,b,c)}.$ 
It follows from the proof of Lemma \ref{lem: omega_type} that each $N_{(\lambda,t,u)}(R_{\omega^\star}) \neq 0,$ and so by the Krull--Schmidt theorem $U(R_{\omega^\star})\neq 0.$ 
By Lemma \ref{lem: brauer_vertex}, every non-projective indecomposable summand of $M_{(2sp,tp,up)}$ therefore has a vertex containing $R_{\omega^\star}.$ 
In the second step of the proof of Theorem \ref{thm: main}, we consider the module $N_{(\lambda,t,u)}(R_{\omega^\star})$ in order to understand $U(R_{\omega^\star}).$\bigskip

\noindent {\bf Second step: The vertices of $N_{(\lambda,t,u)}(R_{\omega^\star}).$}
Recall that we write $K_r$ for the centraliser of $R_r$ in $C_2 \wr S_{rp}.$ 
In this step we show that $N_{(\lambda,t,u)}(R_{\omega^\star})$ is indecomposable as an $FC_{K_r}(R_{\omega^\star})$-module.
It follows that $N_{(\lambda,t,u)}(R_{\omega^\star})$ is an indecomposable $FN_{C_2 \wr S_{rp}}(R_{\omega^\star})$-module, and in Lemma \ref{lem: N_lambda_vertex} we determine its vertex. 
We remark that the group $C_{K_r}(R_{\omega^\star})$ is generated by the set 
\begin{align*}
\{\sigma_i,\tau_i : 1 \le i \le r\} \cup \{\rho_1^{\rho_2\rho_3\ldots\rho_s}\} &\cup \{\rho_i\rho_{i+s} : 1 \le i \le s-1 \mbox{ and } i \neq \lambda_1\}\\ 
&\cup \{\rho_i : 2s+1 \le i < r\},
\end{align*}
and so we have the inclusion $X_{(\lambda,t,u)} \le C_{K_r}(R_{\omega^\star}).$
\begin{lemma}\label{lem: Romega_decomp}
Let $\lambda\in \Lambda(2,s).$ Then the $FN_{K_r}(R_{\omega^\star})$-module $N_{(\lambda,t,u)}(R_{\omega^\star})$ is indecomposable.
\end{lemma}
\begin{proof}
By definition $R_{\omega^\star}$ acts trivially on $N^{\omega^\star}_y,$ and so it follows from \cite[Chapter 8, Corollary 3]{Alperin} that
$$N^{\omega^\star}_y\Res_C \cong F\Ind_{R_{\omega^\star}}^C.$$
This is an indecomposable $FC$-module, and so $N^{\omega^\star}_y$ is an indecomposable $FX_{(\lambda,t,u)}$-module.

Let $(\widetilde{y}(g),\tilde{\gamma},\tilde{\delta}) \in N_{(\lambda,t,u)}(R_{\omega^\star}).$ 
As $C_{K_r}(R_{\omega^\star})$ permutes the $R_{\omega^\star}$-orbits of a fixed size transitively amongst themselves,
it follows from Lemma \ref{lem: tau_action} that there exists
$$\rho \in \langle\rho_1^{\rho_2\rho_3\ldots\rho_s}, \rho_1\rho_{s+1},\ldots,\rho_{s-1}\rho_{2s-1},\rho_{2s+1},\ldots,\rho_{r-1}\rangle$$
such that $\pm(\widetilde{y}(g),\tilde{\gamma},\tilde{\delta})= (y(g^{\rho^{-1}}),\gamma,\delta)\rho,$ where $(y(g),\gamma,\delta)\rho \in N^{\omega^\star}_y.$
Therefore $N^{\omega^\star}_y$ generates $N_{(\lambda,t,u)}(R_{\omega^\star})$ as an $FC_{K_r}(R_{\omega^\star})$-module. 
As there are exactly $\binom{s}{\lambda_1}$
tuples of weight $\lambda$ in $\{-1,1\}^s,$ the second statement of Corollary \ref{rem: tau_special} and Lemma \ref{lem: case_r} imply that the module $N_{(\lambda,t,u)}(R_{\omega^\star})$ has dimension $\binom{s}{\lambda_1} \times \binom{k}{t} \times p^s.$ 
By \cite[Chapter 8, Corollary 3]{Alperin}, we therefore have that
$$N_{{(\lambda,t,u)}}(R_{\omega^\star})\Res_{C_{K_r}(R_{\omega^\star})} \cong N^{\omega^\star}_y\Ind^{C_{K_r}(R_{\omega^\star})}_{X_{(\lambda,t,u)}}.$$
By Lemma \ref{lem: inertial}, the inertial group of $N^{\omega^\star}_y$ in $C_{K_r}(R_{\omega^\star})$ is equal to $X_{(\lambda,t,u)}.$
It follows from \cite[Proposition 3.13.2]{Benson} that $N_{(\lambda,t,u)}(R_{\omega^\star})$ is an indecomposable $FC_{K_r}(R_{\omega^\star})$-module.
\end{proof}

Given $X \subset \{1,2,\ldots,sp\},$ let $C_2 \wr S_X$ be as in \textsection \ref{sec: top_group}.
Also given $x \in \{1,2,\ldots,sp\},$ define $x^* = x + sp.$ 
We remark that this definition of $x^*$ agrees with that of $x^*$ in Example \ref{ex: case_1}, which is the case when $s = 1.$
Given $g \in C_2 \wr S_{\{1,2,\ldots,sp\}},$ let $g^*$ be the permutation in $C_2 \wr S_{\{sp+1,\ldots,2sp\}}$ such that $i^*g^* = (ig)^*.$ 

Also given $\lambda \in \Lambda(2,s),$ we define $J$ to be the group consisting of all elements $gg^*$ such that $g$ is contained in a Sylow $p$-subgroup of $C_2 \wr S_{\{1,\ldots,p\lambda_1\}} \times C_2 \wr S_{\{p\lambda_1+1,\ldots,sp\}}$ with base group $\langle \sigma_1,\ldots,\sigma_s\rangle$. Let $J^+$ be a Sylow $p$-subgroup of 
\[C_2 \wr S_{\{2sp+1,\ldots,(2s+t)p\}} \times C_2 \wr S_{\{(2s+t)p+1,\ldots,rp\}}\]
with base group $\langle \sigma_{2s+1},\ldots,\sigma_r\rangle.$ We define $Q_{(\lambda,t,u)} = J \times J^+.$

By construction, $R_{\omega^\star} \unlhd Q_{(\lambda,t,u)},$ and so $Q_{(\lambda,t,u)} \le N_{C_2 \wr S_{rp}}(R_{\omega^\star}).$
By Lemma \ref{lem: broue_norm} and Lemma \ref{lem: ppermbasis}, there exists a choice of signs $s_v \in \{-1,1\}$ such that 
$\{s_v v: v \in \mathcal B^{R_r}_{(2sp,tp,up)}\}$
is a $p$-permutation basis for $M_{(2sp,tp,up)}(R_{\omega^\star})$ with respect to $Q_{(\lambda,t,u)}.$
Given $v:=(g,\gamma,\delta) \in \mathcal B(\omega^\star),$ let $(h,\tilde{\gamma},\tilde{\delta})$ be a representative for the $Q_{(\lambda,t,u)}$-orbit containing $v.$ 
It follows that for all $E \subseteq \{1,2,\ldots,s\},$ the representative for the $Q_{(\lambda,t,u)}$-orbit containing $(g^{\tau_E},\gamma,\delta)$ can be chosen to be of the form $(h',\tilde{\gamma},\tilde{\delta}).$ 
For distinct summands $w$ and $\widetilde{w}$ of $(y(g),\gamma,\delta),$ it follows that $s_{w} = s_{\widetilde w}.$ 
We can therefore write $s_{(g,\gamma,\delta)}$ in the place of $s_w$ for all such $w,$ and then
\[\{s_{(g,\gamma,\delta)}(y(g),\gamma,\delta) : (g,\gamma,\delta) \in \mathcal B(\omega^\star)\}\]
is a $p$-permutation basis of $N_{(\lambda,t,u)}(R_{\omega^\star})$ with respect to $Q_{(\lambda,t,u)}.$ 
\begin{lemma}\label{lem: N_lambda_vertex}
The module $N_{(\lambda,t,u)}(R_{\omega^\star})$ has a vertex equal to $Q_{(\lambda,t,u)}$.
\end{lemma}
\begin{proof}
Let $y  = y_\lambda.$ 
The element $(f_{sp},\gamma^\star,\delta^\star)$ is a fixed point of $Q_{(\lambda,t,u)}.$ 
As $Q_{(\lambda,t,u)} \le X_{(\lambda,t,u)},$ the element $(y(f_{sp}),\gamma^\star,\delta^\star)$ is also a fixed point of $Q_{(\lambda,t,u)}$. 
Therefore $N_{(\lambda,t,u)}(R_{\omega^\star})$ has a vertex containing $Q_{(\lambda,t,u)}$. 

The element $y(f_{sp})$ is an alternating sum of elements conjugate to $f_{sp}$ in $C_2 \wr S_{rp},$ and so any element in $N_{C_2 \wr S_{rp}}(R_r)$ that fixes $y(f_{sp})$ under the conjugacy action must be contained in $V_{sp}.$  
Indeed, suppose that there exists $h \in Q_{(\lambda,t,u)}$ such that $h \not\in V_{sp}.$
Therefore by definition of $y(f_{sp}),$ it must be the case that $h\tau_E \in V_{sp}$ for some $E \subset \{1,2,\ldots,s\}.$
However $\tau_E$ transposes the $R_{2s}$-orbits
\begin{align*}
\{(j-1)p+1,\ldots,jp\}\\
\{\overline{(j-1)p+1},\ldots,\overline{jp}\},
\end{align*}
for each $j \in E,$ and fixes all other $R_{2s}$-orbits.
As $p$ is odd, it follows that $h$ must act trivially on these orbits. 
The only elements in $N_{C_2 \wr S_{rp}}(R_r)$ that do this are contained in $V_{sp},$ 
which is a contradiction. 

As $Q_{(\lambda,t,u)}$ is the largest $p$-subgroup that is contained in both $X_{(\lambda,t,u)}$ and $V_{sp} \times C_2 \wr S_{tp} \times C_2 \wr S_{up},$ the statement of the lemma now follows by applying Lemma \ref{lem: brauer_vertex}.
\end{proof}

\noindent {\bf Third step: Proof of Theorem \ref{thm: main}.} Given $r \in \N$ such that $rp \le n,$ recall that
$$T'_r = \{(\lambda,t,u) : \lambda \in \Lambda(2,s), 2s+t+u = r \mbox{ and } sp \le a, tp \le b, up \le c\}.$$
We now prove Theorem \ref{thm: main}. We restate the result for the reader's convenience. \bigskip

\noindent {\bf Theorem \ref{thm: main}.} {Let $(a,b,c) \in \N^3_0$ be such that $2a + b + c = n,$ and let $U$ be a non-projective indecomposable summand of $M_{(2a,b,c)}$. 
Then $U$ has a vertex equal to a Sylow $p$-subgroup of 
$$V_{p\lambda} \times C_2\wr S_{tp} \times C_2 \wr S_{up},$$
for some $r \in \N,$ where $rp \le n,$ and $(\lambda,t,u) \in T'_r.$}

\begin{proof}[Proof of Theorem \ref{thm: main}]

Let $r\in\N$ be maximal such that $R_{r}$ is contained in a vertex of $U$.  
By Lemma \ref{lem: tensorfact1}, Lemma \ref{lem: R2r_decomp} and the Krull--Schmidt theorem, there exists $T \subset T'_r,$ and for each $(\lambda,t,u) \in T,$ a summand $W_{(\lambda,t,u)}$ of $M_{(2(a-sp),b-tp,c-up)}$ such that
$$U(R_r) \cong  \bigoplus_{(\lambda,t,u)\in T} N_{(\lambda,t,u)} \boxtimes W_{(\lambda,t,u)},$$
where $s = \size{\lambda}.$

Let $(2s,t,u) \in T$ be such that $s$ is minimal. 
Suppose there exists $(2\tilde{s},\tilde{t},\tilde{u}) \in T^r$ such that $\tilde{s} > s.$
Given $\omega \in \Omega^{(2s; t+u)}$ and $\tilde{\omega} \in \Omega^{(2\tilde{s};\tilde{t}+\tilde{u})},$ the subgroup $R_{\tilde{\omega}}$ cannot contain a conjugate of $R_{\omega}.$
Therefore $N_{(\tilde{\lambda},\tilde{t},\tilde{u})}(R_{\omega}) = 0$ for all $\size{\tilde{\lambda}} = \tilde{s}.$

We therefore now consider $U(R_{\omega}),$ where $\omega \in \Omega^{(2s;t+u)}.$
By Lemma \ref{lem: broue_norm}, there is an isomorphism $U(R_{\omega}) \cong U(R_r)(R_{\omega}),$ and so there exists $S \subseteq T$ such that
$$U(R_{\omega}) \cong \bigoplus_{(\lambda,t,u) \in S} N_{(\lambda,t,u)}(R_{\omega}) \boxtimes W_{(\lambda,t,u)},$$
where $\size{\lambda} = s.$
Let $L = N_{C_2\wr S_{rp}}(R_{\omega}).$
By Lemma \ref{lem: N_lambda_vertex}, each $N_{(\lambda,t,u)}(R_{\omega})$ has a vertex equal to $Q_{(\lambda,t,u)}.$ Let $Q_{(\lambda,t,u)}$ be maximal such that $(\lambda,t,u) \in S.$ 
By Lemma \ref{lem: broue_norm},
\begin{equation}\label{eqn: lambda_decomp}
\begin{array}{ll}
U(Q_{(\lambda,t,u)})&\cong U(R_{\omega})(Q_{(\lambda,t,u)})\\ 
&= \bigoplus_{(\tilde\lambda,\tilde{t},\tilde{u})} N_{(\lambda,t,u)}(R_{\omega})(Q_{(\lambda,t,u)}) \boxtimes W_{(\lambda,t,u)}, 
\end{array}
\end{equation}
where that the sum runs over the $(\tilde\lambda,\tilde{t},\tilde{u}) \in S$ such that $Q_{(\tilde\lambda,\tilde{t},\tilde{u})}$ is a conjugate of $Q_{(\lambda,t,u)}.$ 
Indeed as $N_{(\tilde{\lambda},\tilde{t},\tilde{u})}(R_{\omega})(Q_{(\lambda,t,u)}) \neq 0,$ Lemma \ref{lem: brauer_vertex} says that $Q_{(\lambda,t,u)}$ is contained in a conjugate of $Q_{(\tilde{\lambda},\tilde{t},\tilde{u})},$ say $P.$
If $Q_{(\tilde \lambda, \tilde t, \tilde u)}$ is not a conjugate of $Q_{(\lambda,t,u)},$ then it strictly contains $P.$
However this is a contradiction to the maximality of $Q_{(\lambda,t,u)}.$

As $N_{(\lambda,t,u)}(R_{\omega})(Q_{(\lambda,t,u)}) \neq 0,$ it follows from Lemma \ref{lem: brauer_vertex} that $U$ has a vertex $Q$ containing $Q_{(\lambda,t,u)}.$
Suppose that $Q$ strictly contains $Q_{(\lambda,t,u)}.$ 
As $Q_{(\lambda,t,u)}$ is a $p$-group, there exists $g\in N_Q(Q_{(\lambda,t,u)})$ such that $g\not\in Q_{(\lambda,t,u)}.$ 
The orbits of $Q_{(\lambda,t,u)}$ have length at least $p$ on $\{1,\overline{1},\ldots,rp,\overline{rp}\},$ whereas the orbits of $Q_{(\lambda,t,u)}$ on 
$$\{rp+1,\overline{rp+1},\ldots,n, \overline{n}\}$$
have length 1. As $g$ cannot permute an element in an orbit of length strictly greater than 1 with elements in an orbit of length 1, we can write $g=hh^+,$ where $h\in N_{C_2\wr S_{rp}}(Q_{(\lambda,t,u)})$ and $h^+\in C_2\wr S_{\{rp+1,\ldots,n\}}.$ The only elements in $Q_{(\lambda,t,u)}$ with cycle type either one positive $p$-cycle, or two positive $p$-cycles are those that are contained $R_{\omega}.$ Therefore $N_{C_2\wr S_{rp}}(Q_{(\lambda,t,u)}) \le N_{C_2\wr S_{rp}}(R_{\omega}) = L$, and so $\langle Q_{(\lambda,t,u)}, h \rangle \le N_{L}(Q_{(\lambda,t,u)}).$

Let $\mathcal C$ be a $p$-permutation basis of $N_{(\lambda,t,u)}(R_{\omega})$ with respect to $\langle Q_{(\lambda,t,u)}, g \rangle.$ By Lemma \ref{lem: brauer_vertex}, the group $\langle Q_{(\lambda,t,u)}, g \rangle$ has a fixed point in $\mathcal C.$ It follows from (\ref{eqn: lambda_decomp}) that there exists some $N_{(\lambda,t,u)}(R_{\omega})$ that has a vertex containing $\langle Q_{(\lambda,t,u)}, h\rangle.$ However, we have already seen that each $N_{(\lambda,t,u)}(R_{\omega})$ has vertex equal to $Q_{(\lambda,t,u)}.$ Therefore $h \in Q_{(\lambda,t,u)},$ and so $h^+$ is a non-identity $p$-element of $Q.$ Therefore there exists a power of $h^+$ that is product of positive $p$-cycles with support outside $\{1,\overline{1},\ldots,rp,\overline{rp}\}.$ This contradicts the hypothesis that $r$ is maximal, and so the theorem is proved. 
\end{proof}

\begin{example}\label{ex: main}
In this example, we suppose that $p = 3.$
The module $M_{(54,0,0)}$ is spanned by the conjugates of
\[f_{27}:=(1\ 28)(2\ 29)\ldots(27\ 54)(\overline{1}\ \overline{28})(\overline{2}\ \overline{29})\ldots(\overline{27}\ \overline{54})\]
in $C_2 \wr S_{54}.$ 
In the notation of Theorem 1.1, we have that $r = 9$ and $T'_9 = \Lambda(2,9).$ 
By Theorem 1.1, a non-projective indecomposable summand of $M_{(54,0,0)}$ has a vertex containing a Sylow 3-subgroup of $V_{3\lambda},$ for some $\lambda \in \Lambda(2,9).$ 
In fact we can say more: for every $\lambda \in \Lambda(2,9),$ a Sylow 3-subgroup of $V_{3\lambda}$ contains a conjugate of a Sylow $3$-subgroup of $V_{3(5,4)},$ chosen with the following permutations in its centre:
\[\sigma_1\sigma_{10},\ldots,\sigma_9\sigma_{18}.\]
\end{example}
\section{Decomposition numbers of $C_2 \wr S_n$}\label{sec: decomp}
In this section we prove Theorem \ref{thm: main1}.
In order to do this, we first need to understand how the blocks of $FC_2 \wr S_n$ and the blocks of $FN_{C_2 \wr S_n}(R_r)$ are related. 
We therefore require a description of the blocks of $FN_{C_2 \wr S_n}(R_r),$ which we give in the following section.
\subsection{The blocks of $FN_{C_2 \wr S_n}(R_r)$}
Recall from \eqref{eqn: normaliser_fact} that 
\[N_{C_2 \wr S_n}(R_r) = N_{C_2 \wr S_{rp}}(R_r) \times C_2 \wr S_{\{rp+1,\ldots,n\}}.\]
By Proposition \ref{prop: nakayamconj} the blocks of $FN_{C_2 \wr S_n}(R_r)$ are therefore of the form 
\[b \otimes B((\gamma,v-\widetilde{v}),(\delta,w-\widetilde{w})),\]
where $b$ is a block of $FN_{C_2 \wr S_{rp}}(R_r),$ and $\gamma,\delta$ are $p$-core partitions such that $\size{\gamma} + (v-\widetilde{v})p + \size{\delta}+(w-\widetilde{w})p = n-rp.$ 
It remains to describe the blocks of $FN_{C_2 \wr S_{rp}}(R_r).$ 
\begin{proposition}\label{prop: normaliserblocks}
The blocks of $FN_{C_2 \wr S_{rp}}(R_r)$ are labelled by pairs $(\widetilde{v},\widetilde{w})$ such that $\widetilde v + \widetilde w = r.$ 
\end{proposition}
\begin{proof}
Using the presentation of $N_{C_2 \wr S_{rp}}(R_r)$ given in \textsection \ref{sec: normaliser}, we have that 
\[N_{C_2 \wr S_{rp}}(R_r) \cong C_2^r \rtimes N_{S_{rp}}(R_r).\]
In this case $C_2^r = \langle \tau_1,\ldots,\tau_r\rangle.$ 
Let $\chi_{\widetilde v} \in \Lin(C_2^r)$ be the character such that 
\begin{align*}
\chi_{\widetilde v}(\tau_1) = \cdots = \chi_{\widetilde v}(\tau_{\widetilde{v}}) &= 1\\ 
\chi_{\widetilde v}(\tau_{\widetilde{v}+1}) = \cdots = \chi_{\widetilde v}(\tau_{r}) &= -1.
\end{align*}
Let $\widetilde{w} = r-\widetilde{v}.$ 
Then the stabiliser of $\chi_{\widetilde v}$ in $N_{S_{rp}}(R_r)$ is isomorphic to 
\[H_{(\widetilde v,\widetilde{w})}:=C_{S_{(\widetilde{v},\widetilde{w})p}}(R_r) \rtimes  C_{p-1} ,\] 
and so Theorem \ref{thm: Morita} states that $FN_{C_2 \wr S_{rp}}(R_r)$ and $\bigoplus_{\widetilde{v} = 0}^r FH_{(\widetilde{v},\widetilde{w})}$ are Morita equivalent. The result now follows as $FH_{(\widetilde{v},\widetilde{w})}$ has a unique block by Lemma 2.6 in \cite{BroueFrench} for all $0 \le \widetilde{v} \le r.$  
\end{proof}
We can therefore write $b(\widetilde{v},\widetilde{w})$ for the block of $FN_{C_2 \wr S_{rp}}(R_r)$ labelled by $(\widetilde{v},\widetilde{w}).$ 

Given a $p$-core partition $\gamma = (\gamma_1,\gamma_2,\ldots,\gamma_t)$ and $v \in \N_0,$ define $\gamma+vp$ to be the partition 
$$(\gamma_1+vp,\gamma_2,\ldots,\gamma_t).$$
It is proved in \cite[Lemma 7.1]{wildon2010vertices} that the Specht module $S^{\gamma+vp}$ is a $p$-permutation $FS_{\vert \gamma\vert +vp}$-module. 
If $\delta$ is also a $p$-core partition, then \cite[Proposition 0.2(2)]{broue1985scott} says that $S^{(\gamma+vp,\delta+wp)}$ is a $p$-permutation module for all $v,w \in \N_0.$ 
\begin{proposition}\label{prop: blockC2wrSn}
Fix $v,w \in \N_0.$
Let $\widetilde{v},\widetilde{w} \in \N_0$ be such that $\widetilde{v} \le v, \widetilde{w} \le w$, and $\widetilde{v} + \widetilde{w} = r.$
The $FN_{C_2 \wr S_n}(R_r)$-module $S^{(\gamma+vp,\delta+wp)}(R_r)$ contains a summand lying in the block
$$b(\widetilde{v},\widetilde{w})\otimes B((\gamma,v-\widetilde{v}),(\delta,w-\widetilde{w})).$$
Moreover the blocks $b$ such that $b^{C_2 \wr S_n} = B((\gamma,v),(\delta,w))$ are precisely those of this form.
\end{proposition}

We prove Proposition \ref{prop: blockC2wrSn} by applying Lemma \ref{lem: broue_covering} to $S^{(\gamma+vp,\delta+wp)}$ with respect to $R_r.$ 
Given a $p$-subgroup $Q$ of $C_2 \wr S_n,$ let $U_Q$ denote the kernel of the Brauer morphism from $(S^{(\gamma+vp,\delta+wp)})^{Q}$ to $S^{(\gamma+vp,\delta+wp)}(Q).$
We describe a polytabloid $e_{t_\star}$ that is not contained in $U_Q.$ 
We require the dominance order on standard tableaux, details of which can be found in \cite[Definition 3.11]{James} and \cite[\textsection 3.1]{wildon2010vertices}.

Given a $(\gamma+vp,\delta+wp)$-tableau $t,$ let $\widehat{t}^+$ denote the tableau obtained by replacing each entry $\{x,\overline{x}\}$ in $t^+$ with $x$, and define $\widehat{t}^-$ in the analogous way. 
We write $\widehat{t}$ for the disjoint union of $\widehat{t}^+$ and $\widehat{t}^-.$
Define $t_\star$ to be the tableau such that $\widehat{t_\star}^+$ is the greatest $\gamma+vp$-tableau in the dominance order with entries in $\{1,2,\ldots,\size{\gamma}+vp\},$ and $\widehat{t_\star}^-$ is the greatest $\delta+wp$-tableau in the dominance order with entries in $\{\size{\gamma}+vp+1,\ldots,n\}.$

\begin{lemma}\label{lem: biggest}
Let $Q$ be a $p$-subgroup of $C_2\wr S_n$ with support size $2rp.$
Then the polytabloid $e_{t_\star}$ is not contained in $U_Q.$
\end{lemma}
\begin{proof}
Let $t = t_\star.$ 
Also by definition of the Brauer morphism, we have that $U_Q$ is contained in the subspace $V$ of $S^{(\gamma+vp,\delta+wp)},$ where 
$$V:= \langle e_s + e_sg + \cdots + e_sg^{p-1} : \mbox{$s$ a standard tableau, $g \in Q$}\rangle.$$
We show that $e_t \not\in V.$ 
Suppose, for a contradiction, that $e_t \in V.$ 
Then there exists some $0 \le i \le p-1$ such that $e_t$ has non-zero coefficient in the expression of $e_sg^i$ as a linear combination of standard polytabloids.
Given $g \in Q,$ we can factorise $g = g_+g_-,$ where $g_+ \in C_2 \wr S_{\{1,2,\ldots,\size{\gamma}+vp\}},$ and $g_- \in C_2 \wr S_{\{\size{\gamma}+vp+1,\ldots,n\}}.$ 

Using the bilinearity of the outer tensor product, the polytabloid $e_{t^+}$ has non-zero coefficient in the expression of $e_{s^+}(g_+)^i$ as a linear combination of standard polytabloids.
The analogous statement also holds for $e_{t^-}$ and $e_{s^-}(g_-)^i.$
The action of $Q$ on $e_{t^+}$ (resp. $e_{t^-}$) is equivalent to the action of $\widehat{Q}$ on $e_{\widehat{t}^+}$ (resp. $e_{\widehat{t}^-}).$ 
Therefore it suffices to prove that the polytabloid corresponding to $\widehat{t}^+$ is not contained in the kernel of the Brauer morphism from $(S^{\gamma+vp})^{\widehat{Q}}$ to $S^{\gamma+vp}(\widehat{Q}),$ and that the analogous property holds for the polytabloid corresponding to $\widehat{t}^-.$
This follows from Lemma 5.2 in \cite{wildon2010vertices}.
\end{proof}
Before we prove Proposition \ref{prop: blockC2wrSn}, we introduce one more piece of notation.
Let $M^\lambda$ be the $FS_{\size{\lambda}}$-module corresponding to the action on the cosets of $S_\lambda.$
Given partitions $(\lambda,\mu) \in \mathcal P^2(n),$ we define 
\[M^{(\lambda,\mu)} = (\Inf_{S_{\size{\lambda}}}^{C_2 \wr S_{\size{\lambda}}}M^{\lambda} \boxtimes \widetilde{N}^{\otimes \size{\mu}} \otimes \Inf_{S_{\size{\mu}}}^{C_2 \wr S_{\size{\mu}}} M^\mu)\Ind_{C_2 \wr S_{(\size{\lambda},\size{\mu})}}^{C_2 \wr S_n}.\]
\begin{proof}[Proof of Proposition \ref{prop: blockC2wrSn}]
Let $R_{(\widetilde{v},\widetilde{w})}$ be a conjugate of $R_r$ in the top group, with support such that exactly $\widetilde{v}$ non-trivial orbits of $\widehat{R}_{(\widetilde{v},\widetilde{w})}$ are contained at the end of the first row of $\widehat{t}_\star^+,$ and exactly $\widetilde{w}$ non-trivial orbits of $\widehat{R}_{(\widetilde{v},\widetilde{w})}$ are contained at the end of the first row of $\widehat{t}_\star^-.$

By Lemma \ref{lem: biggest}, the polytabloid $e_{t_\star}$ is not contained in $U_{R_{(\widetilde{v},\widetilde{w})}}.$ 
Therefore the submodule of $S^{(\gamma+vp,\delta + wp)}(R_{(\widetilde{v},\widetilde{w})})$ generated by $e_{t_\star}$, denoted $W$, is non-zero.

Let $s_\star$ be the $(\gamma+(v-\widetilde{v})p,\delta+(w-\widetilde{w})p)$-tableau such that $\widehat{s}_\star^+$ and $\widehat s_\star^-$ are the greatest $\gamma+(v-\widetilde{v})p$-tableau and $\delta + (w - \widetilde{w})p$-tableau in the dominance orders on the tableaux with entries
$$
\begin{array}{l}
\{1,2,\ldots,\size{\gamma}+vp\}\backslash\supp(\widehat{R}_{(\widetilde{v},\widetilde{w})})\\
\{\size{\gamma}+vp+1,\size{\gamma}+vp+2,\ldots,n\}\backslash\supp(\widehat{R}_{(\widetilde{v},\widetilde{w})}),
\end{array}
$$
respectively. 
Let $s$ be the $(\widetilde{v}p,\widetilde{w}p)$-tableau with entries in the row of length $\widetilde{v}p$ agreeing with those at the end of the first row of $t_\star^+,$ and with entries in the row of length $\widetilde{w}p$ agreeing with those at the end of the first row of $t^-_\star.$
The extension of the map $\{s'\}\otimes e_s \mapsto e_t+U$, denoted $\theta,$ is an $F[N_{C_2\wr S_{rp}}(R_{(\widetilde{v},\widetilde{w})})\times C_2 \wr S_{n-rp}]$-module homomorphism from
$$M:= M^{((\widetilde{v}p),(\widetilde wp))}(R_{(\widetilde{v},\widetilde{w})}) \boxtimes S^{(\gamma+(v-\widetilde v)p,\delta + (w - \widetilde w)p)},$$
to $W.$ 
The extension of the map $e_t+U \mapsto \{s'\}\otimes e_s$, denoted $\phi,$ is a well-defined morphism of $F[N_{C_2\wr S_{rp}}(R_{(\widetilde{v},\widetilde{w})})\times C_2 \wr S_{n-rp}]$-modules such that $\phi\theta = \text{id}_M.$ 
Therefore $S^{(\gamma+vp,\delta+wp)}(R_{(\widetilde{v},\widetilde{w})})$ has a submodule isomorphic to $M$.
As $M$ lies in the block $b(\widetilde{v},\widetilde{w}) \otimes B((\gamma,v-\widetilde{v}),(\delta,w-\widetilde{w})),$ there exists a summand of $S^{(\gamma+vp,\delta + wp)}(R_{(\widetilde{v},\widetilde{w})})$ lying in this block, which proves the first statement of the proposition. 
That $B((\gamma,v),(\delta,w))$ corresponds to $b(\widetilde{v},\widetilde{w}) \otimes B((\gamma,v-\widetilde{v}),(\delta,w-\widetilde{w}))$ now follows immediately from Lemma \ref{lem: broue_covering}.

Observe that we have shown if
\[(b(v',w') \otimes B((\gamma',v''),(\delta,w'')))^{C_2 \wr S_n} = B((\gamma,v),(\delta,w)),\]
then $v' + v'' = v$ and $w' + w'' = w.$
In particular $v' \le v$ and $w' \le w.$
Moreover $\gamma' = \gamma$ and $\delta' = \delta.$
This completes the proof of the proposition. 
\end{proof}

\subsection{Proof of Theorem \ref{thm: main1}}
Assume that $M_{(2a,b,c)}$ is defined over the field $\F_p,$ as the results in this section then follow by change of scalars.  
We define $\chi_{(2a,b,c)}$ to be the ordinary character of $M_{(2a,b,c)}$ and $\chi^{(\lambda,\mu)}$ to be the ordinary character of the hyperoctahedral Specht module $S^{(\lambda,\mu)}.$
In the following lemma, we decompose the character $\chi_{(2a,b,c)}$ into its irreducible constituents. 
\begin{lemma}\label{lem: ord_char}
Let $n = 2a + b + c.$
The constituents of the character $\chi_{(2a,b,c)}$ are precisely those $\chi^{(\lambda,\mu)}$ such that the partition $\lambda$ has exactly $b$ odd parts, and the partition $\mu$ has exactly $c$ odd parts, each occuring with multiplicity one.  
\end{lemma}
\begin{proof}
This follows from Propositions 1 and 2 in \cite{baddeley1991models}, and by multiplying through by the ordinary character of the module $\Inf_{S_n}^{C_2 \wr S_n} \sgn_{S_n}.$
\end{proof}

\begin{proposition}\label{prop: proj}
Let $b,c \in \N_0.$
Given $p$-core partitions $\gamma$ and $\delta,$ let $n = \size{\gamma}+w_b(\gamma)p + \size{\delta} + w_c(\delta)p.$
Suppose that if $b,c \ge p,$ then $w_{b-p}(\gamma) \neq w_b(\gamma) -1$ and $w_{c-p}(\delta) \neq w_c(\delta)-1.$ 
Then every summand of $M_{(2a,b,c)}$ lying in the $FC_2 \wr S_n$-block $B((\gamma,w_b(\gamma)),(\delta,w_c(\delta )))$ is projective.
\end{proposition}
\begin{proof}
Suppose for a contradiction that there exists a non-projective indecomposable summand $U$ of $M_{(2a,b,c)}$ in the block $B((\gamma,w_b(\gamma)),(\delta,w_c(\delta)))$. 
By Theorem \ref{thm: main} $U$ has a vertex equal to a subgroup of the form $Q_{(\lambda,t,u)},$
where $\lambda = (\lambda_1,\lambda_2) \vdash s$ and $sp \le a,$ $tp\le b,$ $up\le c.$ 
Let $r = 2s + t + u,$ and so $R_r \le Q_{(\lambda,t,u)}.$
It follows from Lemma \ref{lem: tensorfact1} and Lemma \ref{lem: case_r} that 
\[M_{(2a,b,c)}(R_r) \cong \bigoplus N_{(\lambda,t,u)} \boxtimes M_{(2(a-\size{\lambda}p),b-tp,c-up)},\]
where the sum runs over all $(\lambda,t,u) \in T'_r.$ 
It follows from the Krull--Schmidt theorem that 
\[N_{(\lambda,t,u)} \boxtimes W\, \vert\, U(R_r),\]
for some indecomposable summand $W$ of $M_{(2(a-\size{\lambda}p),b-tp,c-up)}.$
By Lemma \ref{lem: broue_covering}, the block $B((\gamma,w_b(\gamma)),(\delta,w_c(\delta)))$ therefore corresponds to the block containing $N_{(\lambda,t,u)} \boxtimes W.$
The second statement of Proposition \ref{prop: blockC2wrSn} then implies that $W$ lies in a block of the form
$$B:=B((\gamma,w_b(\gamma)-i),(\delta,w_c(\delta)-(r-i))),$$ 
for some $0 \le i \le r.$
By Lemma \ref{lem: ord_char} there exists $S^{(\lambda',\mu')}$ lying in $B$ such that $\lambda'$ has exactly $b-tp$ odd parts, and $\mu'$ has exactly $c-up$ odd parts. 
Adding $tp$ parts of size 1 to $\lambda'$ results in a partition $\lambda$ with $p$-core $\gamma,$ weight $w_b(\gamma)-i+t$ and exactly $b$ odd parts. 
Similarly adding $up$ parts of size 1 to $\mu'$ results in a partition $\mu$ with $p$-core $\delta,$ weight $w_c(\delta)-(r-i)+u$ and exactly $c$ odd parts.
This contradicts the minimality of either $w_b(\gamma)$ or $w_c(\delta),$ unless $(t,u) = (i,r-i).$ 

When $(t,u) = (i,r-i),$ we distinguish two cases.
First suppose that $i \neq 0.$ 
Then adding $(i-1)p$ parts of size 1 to $\lambda'$ results in a partition with $p$-core $\gamma,$ weight $w_b(\gamma)-1$ and $b-p$ odd parts. 
Therefore $w_{b-p}(\gamma) = w_b(\gamma)-1,$ contradicting the hypothesis of the theorem.
In the case that $i = 0,$ we argue in a similar way by adding $(r-1)p$ parts of size 1 to $\mu',$ and contradicting the hypothesis that $w_{c-p}(\delta) \neq w_c(\delta)-1.$ 
\end{proof}
Given $p$-regular partitions $\nu_i$ and $\widetilde{\nu_i},$ let $P^{(\nu_i,\widetilde{\nu_i})}$ denote the projective indecomposable module corresponding to the simple module $D^{(\nu_i,\widetilde{\nu_i})}.$ 
Also let $P^{(\nu_i,\widetilde{\nu_i})}_{\Z_p}$ be the module such that
$$P^{(\nu_i,\widetilde{\nu_i})}_{\Z_p} \otimes_{\Z_p} \F_p = P^{(\nu_i,\widetilde{\nu_i})}.$$
By Brauer reciprocity for projective modules (see for instance \cite[\textsection 9.4]{webb2016course}), the ordinary character of $P^{(\nu_i,\widetilde{\nu_i})}_{\Z_p}$ is
$$\psi^{(\nu_i,\widetilde{\nu_i})} = \sum_{\lambda,\mu} d_{\lambda\nu_i,\mu\widetilde{\nu_i}} \chi^{(\lambda,\mu)},$$
where the decomposition number $d_{\lambda\nu_i,\mu\widetilde{\nu_i}}$ is defined in \textsection \ref{sec: outline}.
It follows from Proposition \ref{prop: nakayamconj} and \cite[Corollary 12.2]{James} that that the sum can be taken over the $(\lambda,\mu) \in \mathcal P^2(n)$ such that $\size{\nu_i} = \size{\lambda}$ with $\lambda \unlhd \nu_i,$ and $\size{\mu} = \size{\widetilde{\nu_i}}$ with $\mu \unlhd \widetilde{\nu_i}.$
\begin{proposition}\label{prop: labels}
Let $b,c \in \N_0.$
Given $p$-core partitions $\gamma$ and $\delta,$ let $n = \size{\gamma}+w_b(\gamma)p + \size{\delta} + w_c(\delta)p.$
Suppose that if $b,c \ge p,$ then $w_{b-p}(\gamma) \neq w_b(\gamma) -1$ and $w_{c-p}(\delta) \neq w_c(\delta)-1.$
Let $\lambda$ and $\mu$ be maximal partitions in $\mathcal{E}_b(\gamma)$ and $\mathcal{E}_c(\delta),$ respectively.
Then $\lambda$ and $\mu$ are both $p$-regular.
\end{proposition}
\begin{proof}
It follows from Proposition \ref{prop: proj} that every summand of the module $M_{(2a,b,c)}$ in the block $B((\gamma,w_b(\gamma)),(\delta,w_c(\delta)))$ is projective. 
Moreover by Lemma \ref{lem: ord_char}, there exists a summand in this block.
Therefore let $$P^{(\nu_1,\widetilde{\nu_1})}_{\F_p},\ldots, P^{(\nu_t,\widetilde{\nu_t})}_{\F_p}$$
be the summands of $M_{(2a,b,c)}$ in the block $B((\gamma,w_b(\gamma)),(\delta,w_c(\delta))).$
Let $M$ denote $M_{(2a,b,c)}$ when defined over $\Z_p.$
It follows from Scott's lifting theorem (see \cite[Corollary 3.11.4]{Benson}) that the summands of $M_{(2a,b,c)}$ can be lifted to summands of $M.$
The ordinary character of the summand of $M_{(2a,b,c)}$ in $B((\gamma,w_b(\gamma))(\delta,w_c(\delta)))$ is equal to $\psi^{(\nu_1,\widetilde{\nu_1})} + \cdots + \psi^{(\nu_t,\widetilde{\nu_t})}.$
It follows from Lemma \ref{lem: ord_char} that
\begin{equation}\label{eqn: max}
\psi^{(\nu_1,\widetilde{\nu_1})} + \cdots + \psi^{(\nu_t,\widetilde{\nu_t})} = \sum_{(\lambda',\mu')}\chi^{(\lambda',\mu')},
\end{equation}
where the sum is over all $(\lambda',\mu') \in \mathcal{E}_b(\gamma) \times \mathcal{E}_c(\delta).$
By Brauer reciprocity, for each $i$ the constituents $\chi^{(\lambda',\mu')}$ of $\psi^{(\nu_i,\widetilde{\nu_i})}$ are such that $\lambda' \unlhd \nu_i$ and $ \mu'\unlhd \widetilde{\nu_i}.$
As $\lambda$ and $\mu$ are maximal, \eqref{eqn: max} gives that $(\nu_i, \widetilde{\nu_i})= (\lambda,\mu)$ for exactly one $i,$ and so the result is proved. 
\end{proof}
Each pair of maximal partitions in $\mathcal E_b(\gamma) \times \mathcal E_c(\delta)$ therefore labels a summand of $M_{(2a,b,c)}$ lying in the block $B((\gamma,w_b(\gamma)),(\delta,w_c(\delta)));$ moreover, every such summand is labelled by a pair of this form.
We now prove Theorem \ref{thm: main1}. 
\begin{proof}[Proof of Theorem \ref{thm: main1}]
Let $P^{(\nu_1,\widetilde{\nu_1})}_{\F_p},\ldots,P^{(\nu_c,\widetilde{\nu_c})}_{\F_p}$ be the summands of $M_{(2a,b,c)}$ lying in the block $B((\gamma,w_b(\gamma)),(\delta,w_c(\delta))),$ all of which are projective. 
It follows from (\ref{eqn: max}) that there exists a set partition $\Lambda_1,\ldots,\Lambda_t$ of $\mathcal{E}_b(\gamma) \times \mathcal{E}_c(\delta)$ such that $(\nu_i,\widetilde{\nu_i}) \in \Lambda_i$ for each $i$ and
$$\psi^{(\nu_i,\widetilde{\nu_i})} = \sum_{(\lambda',\mu')\in \Lambda_i} \chi^{(\lambda',\mu')}.$$
It now follows again using Brauer reciprocity that the column of the decomposition matrix labelled by $(\nu_i,\widetilde{\nu_i})$ has ones in the rows labelled by pairs in $\Lambda_i$ and zeros in all other rows. 
\end{proof}
\section*{Acknowledgements}
This paper was completed by the author under the supervision of Mark Wildon.
The author gratefully acknowledges his support. 
The author also thanks Stephen Donkin and an anonymous referee for their comments on earlier versions of this paper. 

This is a post-peer-review, pre-copyedit version of an article published in Algebras and Representation Theory. 
The final authenticated version is available online at: http://dx.doi.org/10.1007/s10468-018-9839-8.
\bibliographystyle{amsplain}      
\bibliography{bibliography} 

%
%

\end{document}